\documentclass[10pt,letterpaper,twoside,leqno]{article}
\usepackage[T1]{fontenc}
\usepackage{lmodern,microtype}
\usepackage[hcentering,textwidth=5.125in,top=96bp,textheight=594bp,
  headheight=12pt,headsep=12pt,footskip=12pt]{geometry}
\usepackage{amsmath,amssymb,mathtools,amsthm,booktabs,needspace}
\usepackage{titlesec,fancyhdr,float,tabularx,caption,etoolbox}
\numberwithin{equation}{section}

\titleformat{\section}[runin]
  {\normalfont\normalsize\bfseries}{\thesection.}{0.5em}{}[.]
\titlespacing*{\section}{\parindent}{9pt plus 2pt minus 1pt}{0.5em}
\titleformat{\subsection}[runin]
  {\normalfont\normalsize\bfseries}{\thesubsection.}{0.5em}{}[.]
\titlespacing*{\subsection}{\parindent}{6pt plus 2pt minus 1pt}{0.5em}
\titleformat{\subsubsection}[runin]
  {\normalfont\normalsize\bfseries}{\thesubsubsection.}{0.5em}{}[.]
\titlespacing*{\subsubsection}{\parindent}{6pt plus 2pt minus 1pt}{0.5em}
\titleformat{\paragraph}[runin]
  {\normalfont\normalsize\bfseries}{}{0pt}{}[.]
\titlespacing*{\paragraph}{\parindent}{6pt plus 2pt minus 1pt}{0.5em}

\newtheoremstyle{siamplain}{6pt plus 2pt minus 1pt}{6pt plus 2pt minus 1pt}
  {\itshape}{\parindent}{\normalfont\scshape}{.}{0.5em}
  {\thmname{#1}\thmnumber{ #2}\thmnote{\normalfont\ (#3)}}
\theoremstyle{siamplain}
\newtheorem{theorem}{Theorem}[section]
\newtheorem{lemma}[theorem]{Lemma}
\newtheorem{proposition}[theorem]{Proposition}
\newtheorem{corollary}[theorem]{Corollary}
\newtheoremstyle{siamremark}{6pt}{6pt}
  {\normalfont}{\parindent}{\normalfont\scshape}{.}{0.5em}
  {\thmname{#1}\thmnumber{ #2}\thmnote{\normalfont\ (#3)}}
\theoremstyle{siamremark}

\newtheorem{example}[theorem]{Example}
\renewenvironment{proof}[1][\proofname]
  {\par\pushQED{\qed}\normalfont
   \topsep6pt plus 2pt\trivlist\item\relax
   \hspace*{\parindent}\textit{#1.}\enspace\ignorespaces}
  {\popQED\endtrivlist}

\makeatletter
\renewcommand{\maketitle}{%
  \thispagestyle{plain}%
  \begingroup\centering
  \vspace*{10pt}%
  {\fontsize{10}{12}\selectfont\bfseries\MakeUppercase{\@title}\par}%
  \vspace{10pt}%
  {\fontsize{8}{10}\selectfont\MakeUppercase{\@author}\par}%
  \endgroup\vspace{16pt}}
\makeatother
\renewenvironment{abstract}
  {\par\begingroup\fontsize{8}{9.5}\selectfont\setlength{\parindent}{14.4bp}%
   \noindent\hspace*{\parindent}\textbf{Abstract.}\enspace\ignorespaces}
  {\par\endgroup}
\newenvironment{keywords}
  {\par\addvspace{12pt}\begingroup\fontsize{8}{9.5}\selectfont
   \noindent\hspace*{14.4bp}\textbf{Key words.}\enspace\ignorespaces}
  {\par\endgroup}
\newenvironment{AMS}
  {\par\addvspace{12pt}\begingroup\fontsize{8}{9.5}\selectfont
   \noindent\hspace*{14.4bp}\textbf{AMS subject classifications.}\enspace\ignorespaces}
  {\par\endgroup\vspace{6pt}}
\newcommand{\headers}[2]{%
  \fancyhead[CO]{\fontsize{8}{10}\selectfont\MakeUppercase{#1}}%
  \fancyhead[CE]{\fontsize{8}{10}\selectfont\MakeUppercase{#2}}}
\fancypagestyle{plain}{\fancyhf{}\fancyfoot[C]{\footnotesize\thepage}}

\floatstyle{ruled}
\newfloat{algorithm}{tbp}{loa}[section]
\floatname{algorithm}{Algorithm}
\newenvironment{pseudocode}
  {\begin{list}{\arabic{enumi}:}{\usecounter{enumi}
   \setlength{\leftmargin}{18pt}\setlength{\labelsep}{5pt}
   \setlength{\itemsep}{0pt}\setlength{\parsep}{0pt}
   \setlength{\topsep}{2pt}\setlength{\partopsep}{0pt}}}
  {\end{list}}
\usepackage{pgfplots}
\pgfplotsset{compat=1.18}
\usepackage{hyperref}
\hypersetup{hidelinks}
\makeatletter
\renewenvironment{thebibliography}[1]
  {%
    \normalfont\fontsize{8}{9.5}\selectfont
    \titleformat{\section}[block]
      {\normalfont\fontsize{8}{9.5}\selectfont\centering}
      {}{0pt}{\MakeUppercase}
    \titlespacing*{\section}{0pt}{18pt}{10pt}
    \section*{\refname}
    \list{\@biblabel{\@arabic\c@enumiv}}
      {%
        \settowidth{\labelwidth}{\@biblabel{#1}}
        \setlength{\labelsep}{0.6em}
        \setlength{\leftmargin}{\labelwidth}
        \addtolength{\leftmargin}{\labelsep}
        \addtolength{\leftmargin}{1.6em}
        \setlength{\itemindent}{-1.6em}
        \setlength{\itemsep}{0pt}
        \setlength{\parsep}{0pt}
        \setlength{\topsep}{0pt}
        \setlength{\partopsep}{0pt}
        \usecounter{enumiv}
        \let\p@enumiv\@empty
        
      }%
    \sloppy
    \clubpenalty=4000
    \widowpenalty=4000
    \sfcode`\.=\@m
    \urlstyle{same}
  }
  {\endlist}
\makeatother

\newcommand{\epsin}{\varepsilon_{\mathrm{bis}}}
\newcommand{\epsdir}{\varepsilon_{\mathrm{dir}}}
\newcommand{\ip}[2]{\langle #1,#2\rangle}
\newcommand{\nf}[1]{\lVert #1\rVert_F}
\newcommand{\nn}[1]{\lVert #1\rVert_*}
\newcommand{\no}[1]{\lVert #1\rVert_2}
\newcommand{\diag}{\operatorname{diag}}
\newcommand{\rank}{\operatorname{rank}}
\newcommand{\tr}{\operatorname{tr}}
\newcommand{\msign}{\operatorname{msign}}
\newcommand{\R}{\mathbb R}

\newcommand{\cU}{\mathcal U}
\newcommand{\cV}{\mathcal V}
\newcommand{\cF}{\mathcal F}
\newcommand{\cB}{\mathcal B}
\newcommand{\wideD}{\widehat D}
\newcommand{\wideC}{\widehat C}
\newcommand{\kv}{\kappa_v}
\newcommand{\kf}{\kappa_f}
\newcommand{\ph}{\mathrm{ph}}
\newcommand{\sat}{\mathrm{sat}}

\newcommand{\articletitle}{Smoothing Bias and Direction Recovery at Boundary
  Minima of the Nuclear Norm along Matrix Lines}
\title{\articletitle}
\author{Haocheng Wang\textsuperscript{\textasteriskcentered,\textdagger},
  Yuanqi Zhang\textsuperscript{\textasteriskcentered,\textdaggerdbl},
  and Yi Xu\textsuperscript{\textdagger,\textsection}}
\headers{Smoothing Bias along Matrix Lines}{Haocheng Wang, Yuanqi Zhang, and Yi Xu}
\hypersetup{pdftitle={\articletitle},pdfauthor={Haocheng Wang, Yuanqi Zhang, and Yi Xu}}
\date{}
\begin{document}
\maketitle
\begingroup
\renewcommand{\thefootnote}{\fnsymbol{footnote}}
\footnotetext[1]{These authors contributed equally and share first authorship.}
\footnotetext[2]{\raggedright School of Mathematics, Southeast University, Nanjing, China
  \mbox{(Haocheng Wang: \href{mailto:hcwang_28@seu.edu.cn}{\nolinkurl{hcwang_28@seu.edu.cn}}).}}
\footnotetext[3]{\raggedright School of Mathematics and Statistics, Huazhong University of
  Science and Technology, Wuhan, China
  \mbox{(\href{mailto:yuanqizhang375@gmail.com}{\nolinkurl{yuanqizhang375@gmail.com}}).}}
\footnotetext[4]{Corresponding author
  (\href{mailto:yixu@seu.edu.cn}{\nolinkurl{yixu@seu.edu.cn}}).}
\endgroup
\begin{abstract}
We study optimal matrix directions associated with minimizing
$\|G+\lambda\Theta\|_*$, the nuclear norm along a fixed affine line of real
matrices with scalar multiplier $\lambda$. At a rank-deficient minimum,
the compact polar factor can violate the tangency constraint even at the
minimizing multiplier. Spectral smoothing gives a continuous scalar equation
for tangency. We quantify its bias at a unique boundary minimum, where
exactly one one-sided derivative vanishes. An orthogonal block decomposition
shows that pseudo-Huber smoothing with scale $\delta$ shifts the multiplier
by exact order $\delta^{2/3}$ as $\delta\to0$. The Frobenius error relative to
the minimum-Frobenius-norm optimal direction has exact order $\delta^{1/3}$
when some singular values vanish at exactly second order in the multiplier
displacement, and $\delta^{2/3}$ otherwise. A sharp rank bound rules out the
slower rate when $\Theta$ has rank at most two. Huber smoothing gives an
error of exact order $\delta$ for the same limiting direction. Convergent
expansions in powers of $\delta^{1/3}$ also yield a three-scale extrapolation
with $O(\delta)$ error after restoring feasibility.
\end{abstract}
\begin{keywords}
nuclear norm minimization, spectral smoothing, singular value perturbation,
rank deficiency, smoothing bias, direction recovery
\end{keywords}
\begin{AMS}15A18, 15A60, 47A55, 90C25\end{AMS}

\section{Introduction}\label{sec:intro}
A scalar multiplier can be accurately located while a matrix direction
recovered from it remains inaccurate. We investigate this discrepancy for
nuclear-norm minimization along a matrix line, and use the resulting
spectral geometry to improve direction recovery.

\paragraph{The direction problem}
Let $G,\Theta\in\R^{m\times n}$ be given real matrices, with
$\Theta\ne0$. The matrix $G$ specifies a linear objective and $\Theta$ is
the normal to a linear constraint. The unknown is a matrix update direction
$\Phi\in\R^{m\times n}$, used in the additive step $X+\eta\Phi$
from a current matrix $X\in\R^{m\times n}$ with step size $\eta>0$.
For loss minimization we take $G$ to be the negative gradient.
The associated unknown multiplier $\lambda$ is a real scalar. We use the Frobenius inner product
$\ip{A_1}{A_2}=\tr(A_1^TA_2)$ and norm $\nf{A}=\sqrt{\ip{A}{A}}$ for
same-sized matrices $A,A_1,A_2$. Write $q=\min(m,n)$ and denote by
$\sigma_1(A)\ge\cdots\ge\sigma_q(A)\ge0$ the singular values of $A$,
including zeros. The operator norm is $\no{A}=\sigma_1(A)$ and the
nuclear norm is $\nn{A}=\sum_{i=1}^q\sigma_i(A)$.
Nuclear/operator norm duality gives the pair
\begin{equation}\label{eq:duality}
 \max_{\substack{\Phi\in\R^{m\times n}\\
                 \no{\Phi}\le1,\ \ip{\Theta}{\Phi}=0}}\ip{G}{\Phi}
   =\min_{\lambda\in\R}\nn{G+\lambda\Theta}.
\end{equation}
Thus $\lambda$ moves along a known matrix line, while $\Phi$ must satisfy
both a size bound and the tangency equation. We call any maximizing
$\Phi$ an \emph{optimal direction}, or equivalently a dual certificate
for the nuclear-norm minimum.

For linear neural-network layers, maximizing a linear form over a scaled
spectral ball gives the layerwise duality map~\cite[Example~3]{bernstein2025}.
In the spectral sphere optimizer (SSO)~\cite{xie2026}, let
$u_W\in\R^m$ and $v_W\in\R^n$ be unit leading singular vectors of
the weight matrix $W\in\R^{m\times n}$, whose largest singular value
is simple. The sphere has rank-one normal $\Theta=u_Wv_W^T$.
Taking $G$ to be the negative normalized momentum expresses SSO's
unrescaled tangent direction in our additive convention.
At a positive optimal value, the spectral bound is active.
General nonzero normals let us identify the ranks at which slow
recovery first occurs.

\paragraph{Polar directions and the multiplier equation}
SSO~\cite[Section~3.2]{xie2026} uses a multiplier to enforce tangency.
For fixed $\lambda$, adding $\lambda\ip{\Theta}{\Phi}$ to the objective
in \eqref{eq:duality} gives $\ip{G+\lambda\Theta}{\Phi}$.
Over $\no{\Phi}\le1$, one maximizer is the classical compact polar
factor~\cite[Theorem~1.2]{highammehltisseur2010}:
\begin{equation}\label{eq:polar-candidate}
 G+\lambda\Theta=U_\lambda\Sigma_\lambda V_\lambda^T,\qquad
 \msign(G+\lambda\Theta)=U_\lambda V_\lambda^T.
\end{equation}
Here $U_\lambda,V_\lambda$ have respectively $m,n$ rows and orthonormal
columns for the positive singular values in $\Sigma_\lambda$.
The operation $\msign$ replaces those singular values by one and is
zero on the right null space, with $\msign(0)=0$.
This candidate, called the \emph{polar direction}, already satisfies
the norm bound. To impose the remaining tangency constraint, SSO defines
the scalar residual
\begin{equation}\label{eq:polar-residual}
 h(\lambda)=\ip{\Theta}{\msign(G+\lambda\Theta)}
\end{equation}
and searches for a zero by bracketing and bisection.
At a zero, the candidate is feasible and its objective equals the dual
upper bound $\nn{G+\lambda\Theta}$, so it solves \eqref{eq:duality}.

\paragraph{Why multiplier recovery is insufficient}
At a rank-deficient minimizing matrix, a feasible optimal direction may
require a nonzero null-space contribution that the compact polar factor
sets to zero. For example, take
\[
 G=\diag(1,0),\qquad
 \Theta=\tfrac12\begin{pmatrix}1&1\\1&1\end{pmatrix}.
\]
Here $\nn{G+\lambda\Theta}$ equals $\sqrt{1+\lambda^2}$ for $\lambda<0$
and $1+\lambda$ for $\lambda\ge0$, giving the unique minimizing
multiplier $\lambda^\star=0$. Substitution gives
\[
 \msign(G+\lambda^\star\Theta)=\msign(G)=\diag(1,0),
 \qquad h(\lambda^\star)=\tfrac12.
\]
Thus the minimizing multiplier does not make this candidate tangent.
At the same multiplier, every $\diag(1,w)$ with $|w|\le1$ has operator
norm one and maximizes the augmented objective with value one.
Its tangency residual is
$(1+w)/2$, so choosing $w=-1$ gives the optimal direction
$\bar\Phi=\diag(1,-1)$. Recovering a direction therefore requires both
locating the multiplier and selecting its null-space contribution.
Section~\ref{sec:geometry} describes this choice through nuclear-norm
subgradients~\cite{watson1992} before defining the smoothed path.

\paragraph{Smoothing and the boundary case}
At an exact minimizing matrix, its null spaces allow direct completion
of the polar factor, as in the example. Spectral smoothing is an
established alternative~\cite{nesterov2005,beck2012}: it replaces the
polar factor's singular-value coefficients by continuous functions,
so the tangency equation can be solved without a numerical rank threshold.
Pseudo-Huber uses $\sigma/\sqrt{\sigma^2+\delta^2}$, and Huber uses
$\min(\sigma/\delta,1)$ for a singular value $\sigma\ge0$ and scale
$\delta>0$. We determine the direction selected by these approximations
and the error incurred at finite $\delta$.

For a scalar convex objective, a minimum occurs when the interval
between the left and right derivatives contains zero. A \emph{boundary
minimizer} places zero at an endpoint: exactly one derivative vanishes.
It is the limiting case of a strict kink with derivatives of opposite
sign. The $2\times2$ example above has precisely this boundary geometry
and the rank-one normal used in SSO's direction problem.
Here the missing first-order slope makes higher-order spectral changes
control recovery. Some small singular values scale as the multiplier
displacement; others scale as its square and can dominate the direction error.
We call the minimum-Frobenius-norm optimal direction the \emph{center},
denoted by $\bar\Phi$. Same-center recovery means improving the approximation
while retaining this limit.

The following $3\times3$ data illustrate this mechanism:
\begin{equation}\label{eq:slowexample}
 G=\diag(1,0,0),\qquad
 \Theta=\begin{pmatrix}1&0&1\\0&1&0\\1&0&0\end{pmatrix}.
\end{equation}
The exact multiplier is again zero. For pseudo-Huber smoothing with scale
$\delta>0$, the multiplier displacement has order $\delta^{2/3}$,
but the gradient direction approaches $\diag(1,-1,0)$ with error of
order $\delta^{1/3}$ in Frobenius norm. Huber smoothing recovers this
same direction with error of order $\delta$. The computations in
Sections~\ref{sec:ph} and~\ref{sec:profiles} give the corresponding
nonzero coefficients.

\paragraph{Contributions}
For fixed real rectangular matrix lines at nonzero boundary minima:
\begin{enumerate}
\setlength{\itemsep}{0pt}\setlength{\parsep}{0pt}
\item Sections~\ref{sec:structure}--\ref{sec:ph} separate the spectral
blocks governing multiplier displacement and direction error, obtaining
a uniqueness criterion and the two exact pseudo-Huber bias rates.
\item Section~\ref{sec:rank} gives a sharp normal-rank bound.
Rank-one and rank-two normals exclude the slow term; rank-one boundaries
are automatically isolated.
\item Sections~\ref{sec:profiles}--\ref{sec:application} improve recovery
of the center by finite saturation or extrapolation, give a bisection
stopping rule, and bound the error in a constrained update.
\end{enumerate}

\paragraph{Relation to existing analysis}
Affine nuclear-norm minimization and its spectral-ball dual are
classical~\cite{liu2009,recht2010}. Hoheisel and
Paquette~\cite{hoheisel2023} characterize flat nuclear-norm faces and
uniqueness over affine spaces. Our one-dimensional boundary criterion
expresses this geometry through a computable scalar curvature that also
enters the smoothing balance. Second-order variational formulas for the
nuclear norm~\cite{he2026} and spectral differentiation~\cite{lewis2001}
provide local information. Singular-value perturbation theory describes
zero and repeated singular values~\cite{sun1988} and gives a second-order
expansion for the square of the smallest singular value~\cite{stewart1984}.
Here the distinction between exact one-sided
contact and the smoothed balance requires retaining separate spectral
clusters and their matrix coefficients. Orthogonal reductions
preserve the singular values throughout this analysis
\cite{kato1995,bunsegerstner1991}.

Smoothing by a strongly convex dual regularizer is
standard~\cite{nesterov2005,beck2012}. Spectral Huber smoothing is the
Moreau envelope~\cite{parikhboyd2014} of the nuclear norm, with singular-value
thresholding as its proximal map~\cite{cai2010}. Lobos et al.~\cite{lobos2026} develop
Huber-type spectral regularizers and their optimization properties.
We determine the selected direction and its exact bias at a degenerate
boundary, including the change from polynomial tails to finite
saturation. Richardson extrapolation has already been applied to
smoothing parameters~\cite{bach2021}. Our analysis supplies the
fractional-power matrix expansion required here and proves that the
resulting combinations can be made feasible without losing their error
order.

Section~\ref{sec:numerics} tests the predicted rates and the recovery
methods. The appendices give longer proofs and the numerical protocols.

\section{The exact direction and the smoothing target}\label{sec:geometry}
The first task is to specify which optimal direction will be recovered.
We distinguish the given line from the coordinates used to analyze its
minimum, describe all exact directions, and then define the smoothed path.

\subsection{Data, coordinates, and the boundary condition}
Fix the data $G,\Theta$ in \eqref{eq:duality}.
Its scalar objective is coercive because
$\nn{G+\lambda\Theta}\ge|\lambda|\nn{\Theta}-\nn{G}$.
Choose a minimizing multiplier $\lambda^\star\in\R$ and set
\begin{equation}\label{eq:centeredline}
 A^\star=G+\lambda^\star\Theta,\qquad
 s=\lambda-\lambda^\star,\qquad
 \phi(s)=\nn{A^\star+s\Theta}.
\end{equation}
Thus $\phi:\R\to\R$ is the exact centered objective and $s=0$ is a
minimizer by construction. The matrices $G,\Theta$ remain the algorithmic
inputs; the unknown exact solution $\lambda^\star$ and its matrix
$A^\star$ are fixed only for local analysis.

Let $r=\rank A^\star\ge1$. Fix complete orthogonal singular bases
$U_\star=[U_+,U_0]\in\R^{m\times m}$ and
$V_\star=[V_+,V_0]\in\R^{n\times n}$, where
$U_+\in\R^{m\times r}$ and $V_+\in\R^{n\times r}$ span the positive
singular spaces. Set $m_0=m-r$ and $n_0=n-r$; the matrices
$U_0\in\R^{m\times m_0}$ and $V_0\in\R^{n\times n_0}$ span the left
and right null spaces. In these fixed coordinates,
\begin{equation}\label{eq:blocks}
 \begin{aligned}
 U_\star^TA^\star V_\star
   &=\begin{pmatrix}\Sigma&0\\0&0\end{pmatrix},
       \qquad \Sigma=\diag(\sigma_1,\ldots,\sigma_r)>0,\\
 U_\star^T\Theta V_\star
   &=\begin{pmatrix}\Theta_{++}&\Theta_{+0}\\\Theta_{0+}&C\end{pmatrix}.
 \end{aligned}
\end{equation}
Here the positive singular values are ordered
$\sigma_1\ge\cdots\ge\sigma_r>0$, with repetitions allowed.
The block sizes and meanings are
\[
 \begin{array}{c|c|l}
 \text{block}&\text{size}&\text{compression of }\Theta\\ \hline
 \Theta_{++}&r\times r&U_+^T\Theta V_+\\
 \Theta_{+0}&r\times n_0&U_+^T\Theta V_0\\
 \Theta_{0+}&m_0\times r&U_0^T\Theta V_+\\
 C&m_0\times n_0&U_0^T\Theta V_0 .
 \end{array}
\]
We call $(A^\star,\Theta)$ a matrix line.

The block $C$ determines the singular values emerging to first order
from zero. Take its compact singular value decomposition
\begin{equation}\label{eq:C}
 \begin{aligned}
 C&=U_C\Sigma_CV_C^T,\qquad p=\rank C,\\
 \Sigma_C&=\diag(\sigma_{C,1},\ldots,\sigma_{C,p})>0,\\
 P_{C,L}^{\perp}&=I_{m_0}-U_CU_C^T,\qquad
 P_{C,R}^{\perp}=I_{n_0}-V_CV_C^T .
 \end{aligned}
\end{equation}
Thus $U_C\in\R^{m_0\times p}$ and $V_C\in\R^{n_0\times p}$ have
orthonormal columns, $\Sigma_C\in\R^{p\times p}$, and
$\msign(C)=U_CV_C^T$.
The $m_0\times m_0$ and $n_0\times n_0$ orthogonal projectors
$P_{C,L}^{\perp}$ and $P_{C,R}^{\perp}$ have images
$\ker C^T$ and $\ker C$, respectively.
The first-order change of the positive singular-value sum has slope $a$;
new singular values from the null spaces contribute $\beta|s|$ to first order, where
\[
 a=\tr \Theta_{++},\qquad \beta=\nn{C}=\sum_{j=1}^p \sigma_{C,j} .
\]
Thus the null-space contribution has slope $-\beta$ from the left and
$\beta$ from the right. As proved below, the one-sided derivatives of $\phi$ are
$a-\beta$ and $a+\beta$. We study the nonzero endpoint condition
\begin{equation}\label{eq:boundary}
 a=\beta>0.
\end{equation}
It means that the negative-side derivative vanishes while the positive-side
derivative is $2\beta>0$. It implies $p\ge1$ and $r<q$.
The other orientation follows on replacing $(s,\Theta)$ by $(-s,-\Theta)$.
For boundary results set $\sigma_{C,\min}=\min_{1\le j\le p}\sigma_{C,j}>0$.
Residual blocks of size zero are allowed; their norms and ranks are zero.

All matrices in the error statements remain in the original $m\times n$
space. Blocks such as $C$ are coordinates in specified smaller spaces.
Vector norms without subscripts are Euclidean; $I_k$ is the
$k\times k$ identity, $\tr$ is the trace, and $\diag$ forms a block
diagonal matrix. Matrix remainders $O(\cdot)$ use the Frobenius norm.
All asymptotics are for fixed matrix data as the indicated parameter
tends to zero: $O$ denotes an upper bound, $o$ a vanishing ratio,
$\sim$ a ratio tending to one, and $\asymp$ matching upper and lower
bounds by positive constants. The constants and valid neighborhoods
may depend on the fixed data.

\subsection{Complete directions at the boundary}
We use nuclear-norm subgradients to identify the freedom left at a
singular matrix. A matrix $\Phi\in\R^{m\times n}$ is a subgradient at
$A^\star$, written $\Phi\in\partial\nn{\cdot}(A^\star)$, if
\[
 \nn{A}\ge\nn{A^\star}+\ip{\Phi}{A-A^\star}
       \qquad\text{for every }A\in\R^{m\times n}.
\]
In the fixed coordinates of \eqref{eq:blocks}, their set is
\begin{equation}\label{eq:subgradienttool}
 \partial\nn{\cdot}(A^\star)
 =\{U_+V_+^T+U_0WV_0^T:
                 W\in\R^{m_0\times n_0},\ \no{W}\le1\}.
\end{equation}
The positive singular block is fixed, and the contraction $W$ supplies
the null-space completion~\cite{watson1992}. Restriction to the line
pairs these subgradients with $\Theta$.

\begin{proposition}[Optimal directions and their minimum-norm member]\label{prop:cert}
For \eqref{eq:blocks}, $\partial\phi(0)=[a-\beta,a+\beta]$.
Thus zero minimizes $\phi$ if and only if $|a|\le\beta$.
In that case the full optimal direction set in \eqref{eq:duality} is
\begin{equation}\label{eq:facegeneral}
 \cF=\{U_+V_+^T+U_0WV_0^T:
        W\in\R^{m_0\times n_0},\ \no{W}\le1,\ \ip{C}{W}=-a\}.
\end{equation}
Under \eqref{eq:boundary}, the remaining free block acts between
$\ker C$ and $\ker C^T$:
\begin{equation}\label{eq:faceboundary}
 \begin{aligned}
 \cF=\{&U_+V_+^T+U_0(-\msign(C)+W_{\mathrm{res}})V_0^T:\\
 &W_{\mathrm{res}}\in\R^{m_0\times n_0},
                \quad\no{W_{\mathrm{res}}}\le1,\\
 &P_{C,L}^{\perp}W_{\mathrm{res}}P_{C,R}^{\perp}
                         =W_{\mathrm{res}}\}.
 \end{aligned}
\end{equation}
At the boundary, its unique minimum-Frobenius-norm member, the center, is
\begin{equation}\label{eq:center}
 \bar\Phi=U_+V_+^T-U_0\msign(C)V_0^T,\qquad \nf{\bar\Phi}^2=r+p .
\end{equation}
\end{proposition}
\begin{proof}
The affine subgradient rule~\cite[Theorem~23.9]{rockafellar1970} applied to
\eqref{eq:subgradienttool} gives
$\partial\phi(0)=\{a+\ip{C}{W}:\no{W}\le1\}$.
The norm-duality identity used here is
$\max_{\no{W}\le1}\ip{C}{W}=\nn{C}$, attained at $W=\msign(C)$;
$W=-\msign(C)$ attains the negative endpoint. Hence the interval is
$[a-\nn{C},a+\nn{C}]$.
Scalar convex optimality is $0\in\partial\phi(0)$.
A tangent subgradient attains the value $\nn{A^\star}$ in
\eqref{eq:duality}; conversely, an optimal direction satisfies
$\ip{A^\star}{\Phi}=\nn{A^\star}$, so norm duality gives
\[
 \nn{A}\ge\ip{A}{\Phi}
       =\nn{A^\star}+\ip{\Phi}{A-A^\star}
       \qquad(A\in\R^{m\times n}).
\]
This is the subgradient inequality, proving
\eqref{eq:facegeneral} and the duality assertion.

At $a=\beta$, let $u_{C,j}\in\R^{m_0}$ and $v_{C,j}\in\R^{n_0}$ be the
$j$th columns of $U_C,V_C$. The identity
$\sum_j \sigma_{C,j} u_{C,j}^TWv_{C,j}=-\sum_j \sigma_{C,j}$, together with
$u_{C,j}^TWv_{C,j}\ge-1$, gives $u_{C,j}^TWv_{C,j}=-1$ for each $j$.
Consequently
\[
 \|Wv_{C,j}+u_{C,j}\|^2=\|Wv_{C,j}\|^2+1+2u_{C,j}^TWv_{C,j}\le0.
\]
Hence $WV_C=-U_C$; applying the same argument to $W^T$ gives
$U_C^TW=-V_C^T$. Decomposing the domain and range into the active
spaces and their orthogonal complements gives
\[
 W=-\msign(C)+W_{\mathrm{res}},\qquad
 W_{\mathrm{res}}=P_{C,L}^{\perp}WP_{C,R}^{\perp}.
\]
This residual block satisfies the projection invariance in
\eqref{eq:faceboundary}. Conversely, any block satisfying that
invariance acts on domain and range spaces orthogonal to those of
$\msign(C)$. Thus the corresponding $W$ has norm
$\max(1,\no{W_{\mathrm{res}}})$ and obeys $\ip{C}{W}=-\beta$.
This proves both inclusions in \eqref{eq:faceboundary}.
Orthogonality also gives
$\nf{\Phi}^2=r+p+\nf{W_{\mathrm{res}}}^2$.
Its unique minimum occurs at the zero residual block, giving
\eqref{eq:center}.
\end{proof}

\subsection{The pseudo-Huber path and its comparison target}
We begin with pseudo-Huber smoothing because its smooth spectral formula
permits explicit bias expansions. Its gradient approaches the polar
factor on positive singular values while varying continuously through zero.
For any real rectangular input $A$, the approximation and its Frobenius
gradient are
\begin{equation}\label{eq:PHdef}
 \begin{aligned}
 \Psi_\delta(A)&=\sum_i(\sqrt{\sigma_i(A)^2+\delta^2}-\delta),
                     \qquad\delta>0,\\
 Z_\delta(A)&=\nabla\Psi_\delta(A)
           =A(A^TA+\delta^2I)^{-1/2}.
 \end{aligned}
\end{equation}
The sum includes all singular values of the current input, and $I$
has its number of columns. This same convention applies to all smaller
blocks used later. The gradient has the same shape as $A$ and changes a
singular value $\sigma$ into
$\sigma/\sqrt{\sigma^2+\delta^2}$.

Define the direction feasible set
$\cB=\{\Phi\in\R^{m\times n}:\no{\Phi}\le1,\ \ip{\Theta}{\Phi}=0\}$.
Pseudo-Huber smoothing has the Fenchel representation
\begin{equation}\label{eq:PHfenchel}
 \Psi_\delta(A)=\max_{\no{\Phi}\le1}
       \{\ip{A}{\Phi}-\delta\Omega_{\ph}(\Phi)\},\qquad
 \Omega_{\ph}(\Phi)=\sum_i(1-\sqrt{1-\sigma_i(\Phi)^2}).
\end{equation}
Here $\Omega_{\ph}$ is a scalar-valued matrix regularizer on the spectral ball,
extended by $+\infty$ outside it; the subscript denotes pseudo-Huber.
The representation follows from the scalar conjugate identity and
the singular-value conjugacy formula~\cite[Theorem~2.4]{lewis1995}.
It connects the scalar
line minimization to a regularized matrix direction problem.

\begin{proposition}[The selected pseudo-Huber direction]\label{prop:PHselection}
For $\Theta\ne0$ and $\delta>0$, let $s_\delta\in\R$ minimize
$\Psi_\delta(A^\star+s\Theta)$. It exists and is unique.
The associated matrix
\[
 \Phi_\delta=Z_\delta(A^\star+s_\delta\Theta)
\]
is tangent, satisfies $\no{\Phi_\delta}<1$, and is the unique maximizer of
\begin{equation}\label{eq:phprimal}
 \max_{\Phi\in\cB}\{\ip{A^\star}{\Phi}
                         -\delta\Omega_{\ph}(\Phi)\}.
\end{equation}
If \eqref{eq:boundary} holds and $s=0$ is the unique minimizer of $\phi$,
then $s_\delta\to0$ and $\Phi_\delta\to\bar\Phi$ as $\delta\downarrow0$.
\end{proposition}
\paragraph{Selection argument}
The uniform bound
$0\le\nn{A}-\Psi_\delta(A)\le q\delta$ gives existence and convergence
of scalar minimizing points. Strict convexity of the smoothed spectral
trace on any nonconstant matrix line gives uniqueness.
The scalar derivative is $\ip{\Theta}{\Phi_\delta}$, hence vanishes.
Formula~\eqref{eq:PHfenchel} identifies the regularized matrix maximizer.
Every limiting direction minimizes $\Omega_{\ph}$ on $\cF$;
\eqref{eq:faceboundary} makes its unique minimizer the zero residual
completion $\bar\Phi$. Appendix~\ref{app:ph_selection_proof} gives the full
argument, including the spectral Hessian and strong-convexity tools.
The two errors studied below are the scalar displacement $|s_\delta|$
and the matrix bias $\nf{\Phi_\delta-\bar\Phi}$. The reference
$\bar\Phi$ is now fixed even when the exact direction set is nonunique.

\section{The spectral mechanism at a boundary}\label{sec:structure}
The exact direction is now fixed. Its recovery depends on three groups
of singular values: those already positive at $s=0$, those growing as
$|s|\sigma_{C,j}$, and any remaining values of order $s^2$ or smaller.
On the negative side of the boundary, the total first-order change of
the first two groups cancels. Their combined second-order change defines the \emph{partial
curvature} $\kv$; coupling through $\Sigma$ produces a separate
quadratic block $D$ in the residual null spaces. Section~\ref{sec:ph}
uses $\kv$ to locate the smoothed multiplier and $D$ to determine the
slow direction error.

For $1\le i\le r$, write the columns and transposed rows of the two
coupling blocks as
\[
 \theta_{0+,i}=(\Theta_{0+})_{:,i}\in\R^{m_0},\qquad
 \theta_{+0,i}=(\Theta_{+0})_{i,:}^T\in\R^{n_0}.
\]
The residual quadratic block $D$ and the partial curvature $\kv$ are
\begin{align}
 D&=-P_{C,L}^{\perp}\Theta_{0+}\Sigma^{-1}
             \Theta_{+0}P_{C,R}^{\perp}\in\R^{m_0\times n_0},
             \label{eq:D}\\
 \kv&=\sum_{1\le i<j\le r}
       \frac{((\Theta_{++})_{ij}-(\Theta_{++})_{ji})^2}
            {\sigma_i+\sigma_j}\notag\\
 &\quad+\sum_{i=1}^r
   \frac{\|\theta_{0+,i}\|^2+\|\theta_{+0,i}\|^2
      +2\theta_{0+,i}^T\msign(C)\theta_{+0,i}}{\sigma_i}.
      \label{eq:kv}
\end{align}
The block $D$ acts only between the residual spaces of $C$.
To express it in those spaces, fix orthonormal complements
$U_{C,\perp}\in\R^{m_0\times(m_0-p)}$ and
$V_{C,\perp}\in\R^{n_0\times(n_0-p)}$ such that
$[U_C,U_{C,\perp}]$ and $[V_C,V_{C,\perp}]$ are orthogonal.
Its compressed coordinates and the lifts to the original matrix space
are
\begin{equation}\label{eq:residuallifts}
 \begin{aligned}
 D_0&=U_{C,\perp}^TDV_{C,\perp},
       &D&=U_{C,\perp}D_0V_{C,\perp}^T,\\
 \wideD&=U_0DV_0^T,
       &\wideC&=U_0CV_0^T .
 \end{aligned}
\end{equation}
Thus $D_0$ has size $(m_0-p)\times(n_0-p)$, whereas both lifted
matrices have size $m\times n$. Orthogonal compression and lifting
preserve their nonzero singular values. Thus $D_0$, $D$, and $\wideD$
represent the same residual operator in residual, null-space, and original
coordinates; $\wideC$ places $C$ in the original matrix space.

The following reduction makes these roles precise. It uses orthogonal
transformations so that every displayed block represents actual
singular values, including in rectangular and repeated-spectrum cases.

\begin{lemma}[Three spectral scales]\label{lem:normalform}
Let $A(s)=A^\star+s\Theta$ satisfy \eqref{eq:blocks}, with $r,p\ge1$.
There are real-analytic orthogonal matrices
$\cU(s)\in\R^{m\times m}$, $\cV(s)\in\R^{n\times n}$ near $s=0$ such that
\begin{equation}\label{eq:normalform}
 \cU(s)^TA(s)\cV(s)
   =\diag\big(A_+(s),sA_{\mathrm{act}}(s),s^2A_{\mathrm{res}}(s)\big).
\end{equation}
The analytic factors $A_+$, $A_{\mathrm{act}}$, and
$A_{\mathrm{res}}$ have respective sizes $r\times r$, $p\times p$, and
$(m_0-p)\times(n_0-p)$. The active and residual factors are normalized:
their actual blocks carry the powers $s$ and $s^2$ displayed above.
Their initial values are
$A_+(0)=\Sigma$, $A_{\mathrm{act}}(0)=\Sigma_C$, $A_{\mathrm{res}}(0)=D_0$.
The initial bases are
$\cU(0)=[U_+,U_0U_C,U_0U_{C,\perp}]$ and
$\cV(0)=[V_+,V_0V_C,V_0V_{C,\perp}]$.
If $D=0$, the residual block $s^2A_{\mathrm{res}}(s)$ is $O(|s|^3)$.
For the scalar cluster sums
$\ell_+(s)=\nn{A_+(s)}$ and $\ell_{\mathrm{act}}(s)=\nn{A_{\mathrm{act}}(s)}$,
the exact objective on $s<0$ is
\begin{equation}\label{eq:clusterobjective}
 \phi(s)=\ell_+(s)-s\ell_{\mathrm{act}}(s)
                  +s^2\nn{A_{\mathrm{res}}(s)}.
\end{equation}
The derivative of the first two terms, before the residual contribution,
is
\begin{equation}\label{eq:partialbalance}
 \ell_+'(s)-\ell_{\mathrm{act}}(s)-s\ell_{\mathrm{act}}'(s)=a-\beta+\kv s+O(s^2).
\end{equation}
\end{lemma}
\paragraph{Construction}
First separate the positive singular cluster of $A^\star$ from zero.
The remaining block is
$sC-s^2\Theta_{0+}\Sigma^{-1}\Theta_{+0}+O(s^3)$.
Divide by $s$ analytically and separate the positive singular cluster
of $C$; its residual derivative is $D_0$.
The two cluster trace derivatives give \eqref{eq:partialbalance}.
The \hyperlink{proof:normalform}{proof at the start of Appendix~\ref*{app:normalform}}
constructs the moving bases and computes these derivatives.

\subsection{Why the boundary is isolated}
Before calculating a smoothing rate, we need to know whether the exact
multiplier is unique. The next result makes this condition explicit in
the same quantities that govern the bias. Here $\succeq$ denotes the
positive-semidefinite order on real symmetric matrices.

\begin{theorem}[Boundary rigidity]\label{thm:rigidity}
Under \eqref{eq:blocks} and \eqref{eq:boundary},
\begin{equation}\label{eq:rigidityineq}
 \kv\ge2\nn{D}\ge0.
\end{equation}
The exact objective has negative-side contact
\begin{equation}\label{eq:contact}
 \phi(s)=\phi(0)+\tfrac12\kf s^2+O(|s|^3)\quad(s\uparrow0),
 \qquad \kf=\kv+2\nn{D} .
\end{equation}
The full curvature $\kf$ includes the residual contribution
$2\nn{D}$; $\kv$ measures only the first two terms of
\eqref{eq:clusterobjective}. The following are equivalent:
\begin{enumerate}
\item[(i)] $\kv=0$;
\item[(ii)] $\phi$ is constant on $[-\varepsilon,0]$ for some $\varepsilon>0$;
\item[(iii)] zero is a nonunique global minimizer of $\phi$.
\end{enumerate}
In particular, the exact multiplier is unique if and only if $\kv>0$.
\end{theorem}
\begin{proof}
Orthogonal decomposition into the active and residual spaces gives
\begin{equation}\label{eq:squares}
 \begin{aligned}
 &\|\theta_{0+,i}\|^2+\|\theta_{+0,i}\|^2
       +2\theta_{0+,i}^T\msign(C)\theta_{+0,i}\\
 &\quad=\|U_C^T\theta_{0+,i}+V_C^T\theta_{+0,i}\|^2\\
 &\qquad+\|P_{C,L}^{\perp}\theta_{0+,i}\|^2
             +\|P_{C,R}^{\perp}\theta_{+0,i}\|^2.
 \end{aligned}
\end{equation}
Both vectors in the first norm lie in $\R^p$.
Set $X=P_{C,L}^{\perp}\Theta_{0+}\Sigma^{-1/2}\in\R^{m_0\times r}$ and
$Y=\Sigma^{-1/2}\Theta_{+0}P_{C,R}^{\perp}\in\R^{r\times n_0}$, so $D=-XY$.
For each $W\in\R^{m_0\times n_0}$ with $\no{W}\le1$,
$\ip{W}{XY}=\ip{X^TW}{Y}\le\nf{X^TW}\nf{Y}\le\nf{X}\nf{Y}$.
Taking the maximum over such $W$ gives $\nn{XY}\le\nf{X}\nf{Y}$.
Keeping the two residual squares in \eqref{eq:squares}, we obtain
$2\nn{D}\le2\nf{X}\nf{Y}\le\nf{X}^2+\nf{Y}^2\le\kv$.

In \eqref{eq:clusterobjective}, integration of \eqref{eq:partialbalance}, with $a=\beta$, gives
$\ell_+(s)-s\ell_{\mathrm{act}}(s)=\phi(0)+\kv s^2/2+O(|s|^3)$.
Also $\nn{A_{\mathrm{res}}(s)}=\nn{D_0}+O(|s|)=\nn{D}+O(|s|)$.
This proves \eqref{eq:contact}.

If $\kv=0$, every square in \eqref{eq:kv} and \eqref{eq:squares}
vanishes. The residual couplings are zero, and a fixed orthogonal sign
change makes the positive--active block positive definite for small
negative $s$. Its nuclear norm is then its trace,
$\tr\Sigma+s(a-\beta)=\phi(0)$.
Appendix~\ref{app:flat_rigidity} writes this block and verifies positivity
by its Schur complement.

Conversely, $\phi'_+(0)=2\beta>0$ places every other minimizer to the
left of zero. Convexity makes the segment between two minima constant,
so \eqref{eq:contact} gives $\kf=0$.
Together with \eqref{eq:rigidityineq}, this forces $\kv=0$.
\end{proof}

\begin{example}[The quadratic mode in the running example]\label{ex:running}
For \eqref{eq:slowexample}, $A^\star=G$,
$\Sigma=(1)$, $\Theta_{++}=(1)$, $\Theta_{+0}=(0,1)$, $\Theta_{0+}=\Theta_{+0}^T$, and
$C=\diag(1,0)$. Hence $U_C=V_C=(1,0)^T$, $\sigma_{C,\min}=\beta=a=1$,
$D=\diag(0,-1)$, $\kv=2$, and $\kf=4$.
The first and third rows and columns of $A^\star+s\Theta$ form
$\left[\begin{smallmatrix}1+s&s\\s&0\end{smallmatrix}\right]$.
Its smaller eigenvalue is
$((1+s)-\sqrt{(1+s)^2+4s^2})/2=-s^2+O(s^3)$.
Its absolute value is the singular value $s^2+O(|s|^3)$.
The second coordinate supplies the independent singular value $|s|$.
This displays the active and quadratic scales directly.
\end{example}

\section{The pseudo-Huber bias classification}\label{sec:ph}
We can now answer the central quantitative question.
Assume the nonzero boundary \eqref{eq:boundary} is isolated, meaning that
its exact minimizer is unique. Theorem~\ref{thm:rigidity} says exactly
that $\kv>0$.

\paragraph{How the exponents arise}
Write $x=-s>0$ for a test displacement on the negative side.
The positive and active blocks have an unsmoothed derivative
$-\kv x+O(x^2)$ by \eqref{eq:partialbalance}.
The pseudo-Huber coefficient at a positive argument $u$ with $\delta/u\to0$ is
$1-\delta^2/(2u^2)+O(\delta^4/u^4)$.
Applied to the active singular values $x\sigma_{C,j}+O(x^2)$, its total
derivative deficit has leading coefficient
\[
 b_{\ph}=\tfrac12\sum_{j=1}^p \sigma_{C,j}^{-1}>0.
\]
The root therefore balances $\kv x$ against $b_{\ph}\delta^2/x^2$.
This predicts $x$ of order $\delta^{2/3}$.
At that scale the residual block $x^2D_0$ is smaller than $\delta$,
so its gradient is approximately $x^2D_0/\delta$.
A nonzero $D$ consequently produces a direction term of order
$\delta^{1/3}$. The theorem proves this localization, computes the
coefficients, and shows that no additional leading rate occurs when
$D=0$.

Define the positive displacement coefficient $c_{\lambda}$ and its nonnegative
next correction $\nu_{\lambda}$ by
\begin{equation}\label{eq:constants}
 c_{\lambda}=(b_{\ph}/\kv)^{1/3},\qquad
 \nu_{\lambda}=\frac{2c_{\lambda}^3\nf{D}^2}{3\kv}.
\end{equation}
Here $c_\lambda$ multiplies the leading $\delta^{2/3}$ displacement;
$\nu_\lambda$ gives its order-$\delta$ correction from the residual block.
The lifts $\wideD,\wideC$ in \eqref{eq:residuallifts} put all errors in
the original matrix space.

\begin{theorem}[Fixed-instance boundary classification]\label{thm:PH}
Fix a real matrix line satisfying \eqref{eq:boundary} and $\kv>0$,
and let $s_\delta,\Phi_\delta$ be the exact pseudo-Huber path from
Proposition~\ref{prop:PHselection}. There is a fixed direction-bias coefficient
matrix $M_{\Phi}\in\R^{m\times n}$ such that
\begin{align}
 s_\delta&=-c_{\lambda}\delta^{2/3}+\nu_{\lambda}\delta+O(\delta^{4/3}),\label{eq:sexp}\\
 \Phi_\delta&=\bar\Phi+c_{\lambda}^2\wideD\delta^{1/3}
                          +M_{\Phi}\delta^{2/3}+O(\delta),\label{eq:Pexp}\\
 \ip{\wideC}{\Phi_\delta-\bar\Phi}
     &=\kv c_{\lambda}\delta^{2/3}+O(\delta).\label{eq:Cpair}
\end{align}
Both paths have convergent power series in the parameter
$t=\delta^{1/3}$ near zero. The direction bias satisfies the exact laws
\begin{equation}\label{eq:classification}
 \begin{aligned}
 \nf{\Phi_\delta-\bar\Phi}&\sim
 \begin{cases}
  c_{\lambda}^2\nf{D}\,\delta^{1/3},&D\ne0,\\
  \nf{M_{\Phi}}\,\delta^{2/3},&D=0,
 \end{cases}\\
 \nf{M_{\Phi}}&\ge\frac{\kv c_{\lambda}}{\nf{C}}>0\quad(D=0).
 \end{aligned}
\end{equation}
These conclusions include rectangular matrices, empty residual blocks,
and repeated positive singular values. The coefficient $M_\Phi$ combines
motion of the positive and active polar factors, active smoothing, and
the next residual correction; \eqref{eq:M} gives their sum.
\end{theorem}
\begin{proof}
We use the scalar balance from Appendix~\ref{sec:root} and the analytic
block expansions from Appendix~\ref{sec:analytic}. Write
$h_\delta(\lambda^\star+s)=\langle\Theta,Z_\delta(A^\star+s\Theta)\rangle$
for the derivative of the smoothed objective.
At $x=z\delta^{2/3}$, the three-block balance is uniform for $z$
in a compact positive interval:
\[
 h_\delta(\lambda^\star-x)
 =-\kv x+b_{\ph}\delta^2/x^2
     -2\nf{D}^2x^3/\delta+O(\delta^{4/3}).
\]
The sum of the first two terms changes sign at $z=c_\lambda$; the last two are
smaller. Monotonicity brackets the minimizing root between test points
on either side, proving $-s_\delta/\delta^{2/3}\to c_\lambda$.

Set $\delta=t^3$, $s=-t^2z$, and
$\mathcal Z(t,z)=Z_{t^3}(A^\star-t^2z\Theta)$ for $t>0$.
The scaled block formula \eqref{eq:analyticZ} extends this field
analytically through $t=0$, with
$\mathcal Z(0,z)=\bar\Phi$ and $\partial_t\mathcal Z(0,z)=z^2\wideD$.
Both coefficients pair to zero with $\Theta$, since
$\langle\Theta,\bar\Phi\rangle=\langle C,D\rangle=0$.
Thus $h_{\rm sc}(t,z)=t^{-2}\langle\Theta,\mathcal Z(t,z)\rangle$
extends analytically, and the balance above gives
\[
 h_{\rm sc}(t,z)=-\kv z+b_{\ph}z^{-2}
                  -2\nf{D}^2z^3t+O(t^2).
\]
Using $b_{\ph}=\kv c_\lambda^3$, at $(0,c_\lambda)$ we have
\[
 h_{\rm sc}=0,\quad
 \partial_z h_{\rm sc}=-\kv-2b_{\ph}c_\lambda^{-3}=-3\kv<0,
 \quad \partial_t h_{\rm sc}=-2\nf{D}^2c_\lambda^3.
\]
The analytic implicit function theorem gives a unique branch $z(t)$.
The localization identifies it with the minimizing root, and
$z'(0)=-\partial_t h_{\rm sc}/\partial_z h_{\rm sc}=-\nu_\lambda$.
Hence $s_\delta=-t^2z(t)$ gives \eqref{eq:sexp}.

To read the matrix coefficient, retain only the positive and active
polar blocks:
\[
 \Phi_{\rm part}(s)=\cU(s)
 \diag(\msign A_+(s),-\msign A_{\rm act}(s),0)\cV(s)^T,
 \qquad s<0.
\]
This $m\times n$ matrix extends to $\Phi_{\rm part}(0)=\bar\Phi$.
Set $M_{\rm act}=U_0U_C\Sigma_C^{-2}V_C^TV_0^T\in\R^{m\times n}$.
Along $s=-t^2z(t)$, the three contributions are
\[
 \begin{aligned}
 \text{polar motion:}\quad&-c_\lambda t^2\Phi_{\rm part}'(0)+O(t^3),\\
 \text{active smoothing:}\quad&\frac{t^2}{2c_\lambda^2}M_{\rm act}+O(t^3),\\
 \text{residual block:}\quad&(c_\lambda^2t-2c_\lambda\nu_\lambda t^2)\wideD+O(t^3).
 \end{aligned}
\]
Their sum gives \eqref{eq:Pexp} and \eqref{eq:M}. The block expansions also give
$\langle\wideC,\Phi_{\rm part}'(0)\rangle=\langle\wideC,\wideD\rangle=0$
and $\langle\wideC,M_{\rm act}\rangle=2b_{\ph}$.
Thus $\langle\wideC,M_\Phi\rangle=b_{\ph}/c_\lambda^2=\kv c_\lambda>0$,
proving \eqref{eq:Cpair}. For $D=0$, Cauchy--Schwarz yields
$\nf{M_\Phi}\ge\kv c_\lambda/\nf{C}>0$; for $D\ne0$ the leading norm
is $c_\lambda^2\nf{D}>0$. These give both exact rates.
Analyticity of $\mathcal Z(t,z(t))$ gives the convergent matrix series.
\end{proof}
The distinction between scalar and matrix bias is now explicit.
The residual block determines the slow direction term; it enters
the scalar path only through the next coefficient $\nu_{\lambda}$.
The full exact contact curvature $\kf$ cannot replace $\kv$ in $c_{\lambda}$:
they agree precisely when $D=0$.

\begin{example}[The two direction rates]\label{ex:runningPH}
For \eqref{eq:slowexample},
$b_{\ph}=1/2$, $\kv=2$, and $\nf{D}=1$, so
\begin{equation}\label{eq:slowconstants}
 \begin{aligned}
 c_{\lambda}&=4^{-1/3},\qquad \nu_{\lambda}=1/12,\\
 s_\delta&=-4^{-1/3}\delta^{2/3}+\delta/12+O(\delta^{4/3}),\\
 \nf{\Phi_\delta-\bar\Phi}&\sim4^{-2/3}\delta^{1/3}.
 \end{aligned}
\end{equation}
Here $\bar\Phi=\diag(1,-1,0)$.
Appendix~\ref{app:cubic_residual} gives a $4\times4$ example with
$D=0$ but a residual singular value of order $|s|^3$.
Its direction error has order $\delta^{2/3}$ and $\nu_\lambda=0$.
\end{example}

\section{When the normal rank excludes the slow term}\label{sec:rank}
Theorem~\ref{thm:PH} reduces the direction rate to whether $D$ vanishes.
We now give a structural test using the rank of the input normal
$\Theta$. Low rank here refers to $\Theta$, while the returned direction
and the minimizing matrix may have large rank.

\begin{theorem}[Sharp rank protection]\label{thm:rank}
Under \eqref{eq:boundary}, put $d=\rank D$. Then
\begin{equation}\label{eq:rankbound}
 \begin{aligned}
 \rank\Theta
 &\ge p+\rank(P_{C,L}^{\perp}\Theta_{0+})
             +\rank(\Theta_{+0}P_{C,R}^{\perp})\\
 &\ge p+2d .
 \end{aligned}
\end{equation}
The bound is sharp: for every pair of positive integers $(p,d)$,
equality is attained by
a boundary line with $r=d$, $m=n=p+2d$, and $\kv=2d>0$.
Consequently $\rank\Theta\le p+1$ implies $D=0$.
Since the nonzero boundary has $p\ge1$, a normal of rank at most two
satisfies this condition. If the boundary is isolated, its direction
error has precisely order $\delta^{2/3}$ by Theorem~\ref{thm:PH}.
\end{theorem}
\begin{proof}
Use the fixed bases $[U_C,U_{C,\perp}]$ and $[V_C,V_{C,\perp}]$ and define the
active and residual couplings
\[
 \begin{aligned}
 \Theta_{a+}&=U_C^T\Theta_{0+}\in\R^{p\times r},
 &\Theta_{+a}&=\Theta_{+0}V_C\in\R^{r\times p},\\
 \Theta_{\perp+}&=U_{C,\perp}^T\Theta_{0+}
                       \in\R^{(m_0-p)\times r},\\
 \Theta_{+\perp}&=\Theta_{+0}V_{C,\perp}
                       \in\R^{r\times(n_0-p)} .
 \end{aligned}
\]
Reorder the row groups as $(p,r,m_0-p)$ and the column groups as
$(p,r,n_0-p)$. Orthogonal basis changes and these permutations preserve
rank, giving the block matrix
\[
 \begin{pmatrix}
 \Sigma_C&\Theta_{a+}&0\\
 \Theta_{+a}&\Theta_{++}&\Theta_{+\perp}\\
 0&\Theta_{\perp+}&0
 \end{pmatrix}.
\]
Eliminating the invertible $\Sigma_C$ block gives
\begin{equation}\label{eq:rankpivot}
 \rank\Theta=p+\rank
       \begin{pmatrix}\Theta_{++}-\Theta_{+a}\Sigma_C^{-1}\Theta_{a+}&\Theta_{+\perp}\\\Theta_{\perp+}&0\end{pmatrix}.
\end{equation}
Project the column space of the remaining matrix onto its last
$m_0-p$ coordinates. The image has dimension $\rank\Theta_{\perp+}$,
while the kernel contains the $\rank\Theta_{+\perp}$ independent columns
of $\left[\begin{smallmatrix}\Theta_{+\perp}\\0\end{smallmatrix}\right]$.
Rank--nullity bounds the dimension of this column space below by
$\rank\Theta_{\perp+}+\rank\Theta_{+\perp}$.
Since $D_0=-\Theta_{\perp+}\Sigma^{-1}\Theta_{+\perp}$, each rank is at
least $d$.

Appendix~\ref{app:sharp_rank} gives the equality construction.
Finally $p\ge1$ in \eqref{eq:rankbound}, so $\rank\Theta\le2$ forces
$d=0$; Theorem~\ref{thm:PH} gives the rate.
\end{proof}

\begin{corollary}[Rank-one boundaries]\label{cor:rankone}
If $\rank\Theta=1$ and \eqref{eq:boundary} holds, then $p=1$, $D=0$, and
\begin{equation}\label{eq:rankonecurvature}
 \kv\ge\frac{4\beta^2}{\sigma_1}>0,\qquad
 c_{\lambda}\le\frac{\sigma_1^{1/3}}{2\beta}.
\end{equation}
Here $\sigma_1=\sigma_{\max}(A^\star)$ is the largest singular value
of the minimizing matrix.
Both bounds are attained. The multiplier is automatically unique,
and the pseudo-Huber direction error has precisely order $\delta^{2/3}$.
\end{corollary}
\begin{proof}
Factor $\Theta=uv^T$ with $u\in\R^m$, $v\in\R^n$, and write its fixed
singular coordinates as
$u_+=U_+^Tu\in\R^r$, $u_0=U_0^Tu\in\R^{m_0}$, and similarly
$v_+=V_+^Tv\in\R^r$, $v_0=V_0^Tv\in\R^{n_0}$.
The boundary gives $u_+^Tv_+=\|u_0\|\|v_0\|=\beta$.
The active squares in \eqref{eq:squares} therefore satisfy
\[
 \bigl\|\|u_0\|v_++\|v_0\|u_+\bigr\|^2
 =\bigl\|\|u_0\|v_+-\|v_0\|u_+\bigr\|^2+4\beta^2
 \ge4\beta^2.
\]
Their denominators in \eqref{eq:kv} are at most $\sigma_1$, giving
the curvature bound. Then $c_\lambda^3=1/(2\beta\kv)$ gives the
coefficient bound. Appendix~\ref{app:sharp_rank} supplies the block
identities and examples attaining both bounds.
\end{proof}
This is the return from the general classification to the SSO
direction subproblem: its rank-one normal excludes the quadratic slow
term at every nonzero boundary. The two recovery methods in the next
section improve the remaining $\delta^{2/3}$ bias while keeping
the selected optimal direction.

\section{Improving recovery without changing the limiting direction}
\label{sec:profiles}\label{sec:algorithms}
Theorem~\ref{thm:PH} identifies the slow term and the direction it approaches.
We now use this information in two ways: change the scalar smoothing
profile so that its derivative reaches one at a finite argument, or
cancel terms of the original pseudo-Huber path by extrapolation.
Both constructions approach the same $\bar\Phi$.
We first establish this common target, then give the improvements and
their finite-accuracy implementation.

\subsection{Spectral profiles and their common target}
A smoothing profile specifies how each singular value contributes to
the recovered direction: its derivative $f$ sends a singular value
$\sigma$ to the coefficient $f(\sigma/\delta)$.
Formally, a \emph{profile} is an even convex $C^1$ function
$\psi:\R\to\R$ with $\psi(0)=0$ and
$f=\psi':[0,\infty)\to[0,1]$.
Pseudo-Huber uses $\psi_{\ph}(z)=\sqrt{1+z^2}-1$.
We call $\psi$ \emph{admissible} when $f(0)=0$,
$f$ increases strictly until it reaches one and is one thereafter
(the threshold may be infinite),
$\lim_{z\to\infty}(z-\psi(z))<\infty$, and its Fenchel conjugate
\[
 \psi^*(u)=\sup_{z\in\R}\{uz-\psi(z)\}
\]
is finite and $\mu$-strongly convex on $[-1,1]$ for some $\mu>0$,
and is infinite outside this interval.
Here strong convexity means that
$u\mapsto\psi^*(u)-\mu u^2/2$ is convex on $[-1,1]$.

For $A,\Phi\in\R^{m\times n}$, define the scalar-valued spectral functions
\begin{equation}\label{eq:generalprofile}
 \Psi_\delta^\psi(A)=\delta\sum_{i=1}^q\psi(\sigma_i(A)/\delta),
 \qquad
 \Omega_\psi(\Phi)=
 \begin{cases}
 \displaystyle\sum_{i=1}^q\psi^*(\sigma_i(\Phi)),&\no{\Phi}\le1,\\
 +\infty,&\no{\Phi}>1 .
 \end{cases}
\end{equation}
The gradient of $\Psi_\delta^\psi$ is explicit: if
$A=U\diag(\sigma_1(A),\ldots,\sigma_q(A))V^T$ is an SVD with
$U\in\R^{m\times q}$ and $V\in\R^{n\times q}$, including zero
singular values, then
\[
 Z_\delta^\psi(A):=\nabla\Psi_\delta^\psi(A)
 =U\diag\big(f(\sigma_1(A)/\delta),\ldots,
               f(\sigma_q(A)/\delta)\big)V^T .
\]
The gradient uses the Frobenius inner product; $f(0)=0$ makes the
expression independent of the zero-singular-value bases.
The same definitions apply to each rectangular block, with $q$ and
the identity matrices determined by that block's dimensions.

Write $s_\delta^\psi$ for any minimizer of
$s\mapsto\Psi_\delta^\psi(A^\star+s\Theta)$, and set
$\Phi_\delta^\psi=Z_\delta^\psi(A^\star+s_\delta^\psi\Theta)$.
A missing profile superscript denotes pseudo-Huber, for both exact
and computed directions; its regularizer is $\Omega_{\ph}$ in
\eqref{eq:PHfenchel}. The next proposition identifies the same boundary
center even when the optimal direction is nonunique.

\begin{proposition}[Common boundary selection]\label{prop:selection}
For an admissible $\psi$ and $\Theta\ne0$, a minimizing
$s_\delta^\psi$ exists for every $\delta>0$.
All such minimizers give the same $\Phi_\delta^\psi$, the unique
solution of
\begin{equation}\label{eq:regularizedprimal}
 \max_{\Phi\in\cB}
 \{\ip{A^\star}{\Phi}-\delta\Omega_\psi(\Phi)\},
 \qquad
 \cB=\{\Phi\in\R^{m\times n}:\no{\Phi}\le1,\
                              \ip{\Theta}{\Phi}=0\}.
\end{equation}
If \eqref{eq:boundary} holds and $\kv>0$, then every
$s_\delta^\psi\to0$ and $\Phi_\delta^\psi\to\bar\Phi$ as
$\delta\downarrow0$.
\end{proposition}

The spectral Fenchel formula turns the scalar minimizing condition
into tangency. On the exact boundary face, every admissible
$\Omega_\psi$ is minimized by setting its free residual block to zero.
Appendix~\ref{app:selection} proves these two facts and the convergence.

\subsection{Changing the tail or imposing finite saturation}
The active singular values sample $f$ at arguments of size
$|s|/\delta$, while the quadratic residual samples it at size
$s^2/\delta$. This separates the two profile properties relevant to
the direction error: the tail deficit $1-f(z)$ as $z\to\infty$, and
the slope of $f$ at zero.

\begin{theorem}[Polynomial-tail recovery]\label{thm:tail}
Assume \eqref{eq:boundary}, $\kv>0$, and an admissible profile
with tail exponent $\rho>1$, coefficient $\gamma>0$, and slope
$\alpha_0>0$ at zero, satisfying
\begin{equation}\label{eq:tailhyp}
 \begin{aligned}
 1-f(z)&=\gamma z^{-\rho}+o(z^{-\rho})&& (z\to\infty),\\
 f(z)&=\alpha_0 z+o(z)&& (z\downarrow0).
 \end{aligned}
\end{equation}
Every minimizing multiplier satisfies
$|s_\delta^\psi|\asymp\delta^{\rho/(\rho+1)}$, and
\begin{equation}\label{eq:tailrates}
 \nf{\Phi_\delta^\psi-\bar\Phi}\asymp
 \begin{cases}
 \delta^{(\rho-1)/(\rho+1)},&D\ne0,\\
 \delta^{\rho/(\rho+1)},&D=0 .
 \end{cases}
\end{equation}
\end{theorem}

The active deficit changes the scalar balance to
$\kv x\asymp\delta^\rho/x^\rho$, where $x=-s_\delta^\psi>0$;
the residual direction still has size $x^2/\delta$.
Appendix~\ref{app:tail} gives the matching matrix lower bounds.
For example, taking
$f_\rho(z)=z/(1+z^\rho)^{1/\rho}$ for $z\ge0$ and
$\psi_\rho(z)=\int_0^{|z|}f_\rho(u)\,du$
gives $\gamma=1/\rho$, $\alpha_0=1$, and $\mu=1$:
$f_\rho'(z)=(1+z^\rho)^{-1-1/\rho}\le1$ makes its
conjugate 1-strongly convex.
Pseudo-Huber is $\rho=2$.

A finite saturation threshold changes the displacement to order
$\delta$: the weakest active singular value reaches the constant part
of $f$ when $|s|\sigma_{C,\min}/\delta$ reaches that threshold.

\begin{theorem}[Linear recovery by finite saturation]\label{thm:saturation}
Assume \eqref{eq:boundary}, $\kv>0$, and an admissible profile such that,
for a fixed $z_\sat>0$,
$f(z)<1$ for $0\le z<z_\sat$ and $f(z)=1$ for $z\ge z_\sat$.
With $\sigma_{C,\min}=\min_{1\le j\le p}\sigma_{C,j}$, every minimizing multiplier obeys
\begin{align}
 s_\delta^\psi&=-\frac{z_\sat}{\sigma_{C,\min}}\delta+o(\delta),
 &\nf{\Phi_\delta^\psi-\bar\Phi}&\asymp\delta,\label{eq:satmain}\\
 \frac{\ip{\wideC}{\Phi_\delta^\psi-\bar\Phi}}{\delta}
 &\longrightarrow\frac{\kv z_\sat}{\sigma_{C,\min}}>0 .\label{eq:satpair}
\end{align}
For Huber's profile
\[
 \psi_{\rm hub}(z)=
 \begin{cases}z^2/2,&|z|\le1,\\ |z|-1/2,&|z|\ge1,\end{cases}
 \qquad f_{\rm hub}(z)=\min(z,1)\quad(z\ge0),
\]
the multiplier is unique for sufficiently small $\delta$, and
$s_\delta^{\rm hub}=-\delta/\sigma_{C,\min}+O(\delta^2)$.
Here the profile superscript $\mathrm{hub}$ abbreviates $\psi_{\rm hub}$.
\end{theorem}

\begin{proof}
Put $\tau_0=z_\sat/\sigma_{C,\min}$.
At $s=-\tau\delta$ with fixed $0<\tau<\tau_0$, the derivative of
$\Psi_\delta^\psi(A^\star+s\Theta)$ tends to
$\beta-\sum_j\sigma_{C,j}f(\tau\sigma_{C,j})>0$.
For $\tau>\tau_0$ it equals $-\kv\tau\delta+O(\delta^2)<0$.
Monotonicity therefore places every root at
$x=-s_\delta^\psi=\tau_0\delta+o(\delta)$.

Write the active singular values as $xk_j(-x)$, where
$k_j(0)=\sigma_{C,j}$, and define their nonnegative coefficient deficits
$d_j=1-f(xk_j(-x)/\delta)$ and $d_{\max}=\max_jd_j$.
The derivative decomposition gives
\[
 \sum_j\sigma_{C,j}d_j=\kv x+O(x^2+xd_{\max}).
\]
For a fixed remainder bound $K(x^2+xd_{\max})$ with $K\ge0$,
$(\sigma_{C,\min}-Kx)d_{\max}\le\kv x+Kx^2$.
Taking $Kx\le\sigma_{C,\min}/2$ gives $d_{\max}=O(x)$.
Frame motion, active deficits, and the residual direction
$O(x^2/\delta)$ are consequently all $O(\delta)$.
The fixed $\wideC$ pairing removes the first-order frame motion and gives
$\ip{\wideC}{\Phi_\delta^\psi-\bar\Phi}=\kv x+O(\delta^2)$,
proving \eqref{eq:satpair} and, by Cauchy--Schwarz, the nonzero lower bound.
For Huber, an active mode stays unsaturated at a root; its deficit
$1-xk_j(-x)/\delta=O(\delta)$ gives the stated multiplier correction.
Appendix~\ref{app:saturation} proves the block estimates and uniqueness
across the clipping points.
\end{proof}
For the running example \eqref{eq:slowexample}, the direction error is
$\sqrt7\,\delta(1+o(1))$, compared with $4^{-2/3}\delta^{1/3}$ for
pseudo-Huber. Appendix~\ref{app:saturation} gives its matrix expansion;
Appendix~\ref{app:edge} treats higher-order approach to a saturation edge.

\subsection{Extrapolating the original path}\label{sec:extrap}
The second route keeps the pseudo-Huber profile.
Theorem~\ref{thm:PH} supplies a convergent matrix series in
$t=\delta^{1/3}$, so Richardson extrapolation
\cite{bach2021,sidi2003} can cancel its first terms.
Because the weights include negative entries, we restore spectral
feasibility by the radial map
$\mathcal N(Q)=Q/\max(1,\no{Q})$ for $Q\in\R^{m\times n}$.

\begin{corollary}[Extrapolation with radial repair]\label{cor:extrap}
Under Theorem~\ref{thm:PH}, define the raw two- and three-scale matrices
\begin{equation}\label{eq:extrapdef}
 Q_\delta^{(2)}=2\Phi_{\delta/8}-\Phi_\delta,\qquad
 Q_\delta^{(3)}
   =(\Phi_\delta-6\Phi_{\delta/8}+8\Phi_{\delta/64})/3 .
\end{equation}
The superscript counts the scales used.
Both $\mathcal N(Q_\delta^{(2)})$ and
$\mathcal N(Q_\delta^{(3)})$ belong to $\cB$, and
\begin{equation}\label{eq:extraprates}
 \nf{\mathcal N(Q_\delta^{(2)})-\bar\Phi}=O(\delta^{2/3}),
 \qquad
 \nf{\mathcal N(Q_\delta^{(3)})-\bar\Phi}=O(\delta).
\end{equation}
\end{corollary}

At the scales $\delta/8$ and $\delta/64$, the series parameter is
$t/2$ and $t/4$.
The weights preserve the constant term since $-1+2=1$ and
$(1-6+8)/3=1$. Their cancelled coefficients are
\[
 -1+2/2=0,\qquad
 (1-6/2+8/4)/3=0,\qquad (1-6/4+8/16)/3=0.
\]
Thus the raw errors are $O(t^2)$ and $O(t^3)$, respectively.
Spectral feasibility of $\bar\Phi$ gives
$\max(1,\no{Q})-1\le\nf{Q-\bar\Phi}$ and hence
\begin{equation}\label{eq:repairbound}
 \nf{\mathcal N(Q)-\bar\Phi}
 \le(1+\nf{\bar\Phi})\nf{Q-\bar\Phi}.
\end{equation}
Thus repair preserves the improved order and, being radial, preserves
tangency. Appendix~\ref{app:extrapolation} supplies the complete weight
and repair argument, including the general Vandermonde construction.

\subsection{Computing a direction by bisection}
We now distinguish the exact smoothed direction $\Phi_\delta^\psi$
from a computed approximation $\widetilde\Phi_\delta^\psi$.
In particular, $\widetilde\Phi_\delta$ denotes the computed
pseudo-Huber direction.
The input matrices are $G,\Theta\in\R^{m\times n}$, not the unknown
centered matrix $A^\star$.
A smoothed minimizing multiplier in the input coordinate is
$\lambda_\delta^\psi=\lambda^\star+s_\delta^\psi$, since
\begin{equation}\label{eq:R3coordinates}
 G+\lambda\Theta=A^\star+s\Theta,\qquad s=\lambda-\lambda^\star.
\end{equation}
The scalar derivative used by the computation is the tangency residual
\begin{equation}\label{eq:R2oracle}
 h_\delta^\psi(\lambda)=
 \ip{\Theta}{Z_\delta^\psi(G+\lambda\Theta)},\qquad \lambda\in\R .
\end{equation}
Strong convexity of $\Omega_\psi$ yields the matrix gradient bound
\begin{equation}\label{eq:R2spectralLip}
 \nf{Z_\delta^\psi(A_1)-Z_\delta^\psi(A_2)}
 \le\frac{\nf{A_1-A_2}}{\mu\delta}
 \qquad(A_1,A_2\in\R^{m\times n}).
\end{equation}
Indeed, the Fenchel maximizers satisfy
$A_i\in\delta\partial\Omega_\psi(Z_\delta^\psi(A_i))$.
The strong-monotonicity inequality for these two subgradients,
followed by Cauchy--Schwarz, gives \eqref{eq:R2spectralLip};
Appendix~\ref{app:oracle} writes out the inequality.
Pseudo-Huber and Huber both have $\mu=1$.

With exact evaluations, stopping bisection early leaves an error relative
to $\Phi_\delta^\psi$. To bound it by $\epsin>0$, Algorithm~\ref{alg:oracle}
uses the interval-width threshold $2\delta\epsin/\nf{\Theta}$ from
\eqref{eq:R2spectralLip}.

\begin{algorithm}[htbp]
\caption{Smoothed direction with matrix tolerance $\epsin$}
\label{alg:oracle}
\begin{pseudocode}
\item \textbf{Input:} $G,\Theta\in\R^{m\times n}$, $\Theta\ne0$,
 $\delta>0$, pseudo-Huber or Huber profile, and $\epsin>0$.
\item Set $R_{\rm br}=1$; double it while
 $h_\delta^\psi(-R_{\rm br})>0$ or $h_\delta^\psi(R_{\rm br})<0$.
 Set $\lambda_-=-R_{\rm br}$ and $\lambda_+=R_{\rm br}$.
\item \textbf{While} $\lambda_+-\lambda_->2\delta\epsin/\nf{\Theta}$:
\item \hspace*{1em}Set $\lambda_{\rm mid}=(\lambda_-+\lambda_+)/2$ and
 $h_{\rm mid}=h_\delta^\psi(\lambda_{\rm mid})$.
\item \hspace*{1em}If $h_{\rm mid}=0$, return
 $Z_\delta^\psi(G+\lambda_{\rm mid}\Theta)$.
\item \hspace*{1em}If $h_{\rm mid}<0$, set $\lambda_-=\lambda_{\rm mid}$;
 otherwise set $\lambda_+=\lambda_{\rm mid}$.
\item \textbf{End while.} Return
 $\widetilde\Phi_\delta^\psi
 =Z_\delta^\psi(G+\tfrac12(\lambda_-+\lambda_+)\Theta)$.
\end{pseudocode}
\end{algorithm}

\begin{proposition}[Matrix accuracy from bisection]\label{prop:R2oracle}
With exact scalar and matrix evaluations,
Algorithm~\ref{alg:oracle} terminates and returns
$\nf{\widetilde\Phi_\delta^\psi-\Phi_\delta^\psi}\le\epsin$.
It also gives
$|\ip{\Theta}{\widetilde\Phi_\delta^\psi}|
 \le\nf{\Theta}\epsin$ and $\no{\widetilde\Phi_\delta^\psi}\le1$.
Starting from a bracket of width $w_{\rm br}$, at most
$\max\{0,\lceil\log_2(w_{\rm br}\nf{\Theta}/(2\delta\epsin))\rceil\}$
bisections are needed.
\end{proposition}
Appendix~\ref{app:oracle} proves the bracket and error bounds.

\paragraph{Returning an exactly feasible direction}
A stopped bisection may leave a nonzero tangency residual.
For a raw matrix $Q\in\R^{m\times n}$, first project in the Frobenius
matrix space onto the hyperplane normal to $\Theta$, then normalize:
\begin{equation}\label{eq:R2repair}
 \Pi_\Theta(Q)=Q-\frac{\ip{\Theta}{Q}}{\nf{\Theta}^2}\Theta,\qquad
 \mathcal R_\Theta(Q)=\mathcal N(\Pi_\Theta(Q)).
\end{equation}
Projection restores tangency and normalization enforces the spectral
bound. Both maps fix $\bar\Phi$.

\begin{algorithm}[t]
\caption{Same-center recovery with feasibility repair}
\label{alg:recovery}
\begin{pseudocode}
\item \textbf{Input:} $G,\Theta\in\R^{m\times n}$, $\Theta\ne0$,
 $\delta>0$, $\epsin>0$, and mode H (Huber) or E (extrapolation).
\item In mode H, run Algorithm~\ref{alg:oracle} with Huber at
 $\delta$ and tolerance $\epsin$; set
 $Q=\widetilde\Phi_\delta^{\rm hub}$.
\item In mode E, run Algorithm~\ref{alg:oracle} with pseudo-Huber at
 $\delta,\delta/8,\delta/64$, each with tolerance $\epsin$; set
 $Q=(\widetilde\Phi_\delta-6\widetilde\Phi_{\delta/8}
        +8\widetilde\Phi_{\delta/64})/3$.
\item \textbf{Output:} $\widehat\Phi=\mathcal R_\Theta(Q)$.
\end{pseudocode}
\end{algorithm}

\begin{proposition}[Total recovery error]\label{prop:R2recovery}
Assume exact scalar and matrix computations and a unique minimizing
$\lambda^\star$ whose matrix satisfies \eqref{eq:boundary}.
There exist $C_{\rm rec},\delta_0>0$ depending only on the fixed input
line and recovery mode. For every $0<\delta<\delta_0$ and $\epsin>0$,
Algorithm~\ref{alg:recovery} returns
$\ip{\Theta}{\widehat\Phi}=0$, $\no{\widehat\Phi}\le1$, and
\begin{equation}\label{eq:R2finiteaccuracy}
 \nf{\widehat\Phi-\bar\Phi}
 \le(1+\sqrt q)(C_{\rm rec}\delta+\omega_{\rm mode}\epsin),
 \qquad \omega_{\rm H}=1,\quad\omega_{\rm E}=5 .
\end{equation}
In particular $\epsin=O(\delta)$ gives total direction error $O(\delta)$.
For any raw matrices $Q,Q'\in\R^{m\times n}$,
\begin{equation}\label{eq:R2repairLip}
 \nf{\mathcal R_\Theta(Q)-\mathcal R_\Theta(Q')}
 \le(1+\sqrt q)\nf{Q-Q'} .
\end{equation}
\end{proposition}
The sum of absolute extrapolation weights is $\omega_{\rm E}=(1+6+8)/3=5$.
It bounds the amplification of the three bisection errors.
Appendix~\ref{app:recovery} combines these errors with the exact-path
bias and proves the repair bound.

\section{Application to a constrained matrix update}\label{sec:application}
The recovery results return a feasible matrix direction with a controlled
error. We now show which quantities in one constrained update inherit
that control, and why the same boundary center is a useful reference.

Let $\mathcal L$ be a real-valued differentiable function on a
neighborhood in $\R^{m\times n}$. At a current feasible matrix $X$ for
the fixed linear constraint $\ip{\Theta}{X}=b$, where $b\in\R$ is given,
put $G=-\nabla\mathcal L(X)$. The gradient is taken in the Frobenius
inner product. This convention uses the additive update
$X^+=X+\eta\Phi$ with step length $\eta\ge0$.
At this step the optimal linear decrease is
\begin{equation}\label{eq:R2directionproblem}
 v^\star=\max_{\Phi\in\cB}\ip{G}{\Phi}
             =\min_{\lambda\in\R}\nn{G+\lambda\Theta}.
\end{equation}
Tangency makes the additive update exactly feasible for the linear
constraint. In particular, a scalar measurement $u^TXv=b$ has
$\Theta=uv^T$ for given nonzero $u\in\R^m$, $v\in\R^n$.
Its rank-one normal permits the direct use of
Corollary~\ref{cor:rankone} at a nonzero boundary.

\begin{proposition}[Direction error and one-step bounds]\label{prop:R2application}
Assume the line in \eqref{eq:R2directionproblem} satisfies
Proposition~\ref{prop:R2recovery}. For a feasible computed direction
$\widehat\Phi\in\cB$, define its actual error
$e_{\rm dir}=\nf{\widehat\Phi-\bar\Phi}\ge0$.
Then for every $\eta\ge0$,
\begin{equation}\label{eq:R2applicationbound}
 0\le v^\star-\ip{G}{\widehat\Phi}\le\nf{G}\,e_{\rm dir},\qquad
 \nf{(X+\eta\widehat\Phi)-(X+\eta\bar\Phi)}=\eta e_{\rm dir}.
\end{equation}
Suppose $\nabla\mathcal L$ is $L_\nabla$-Lipschitz in Frobenius norm
on a convex neighborhood containing the segment
$\{X+t\widehat\Phi:0\le t\le\eta\}$, with $L_\nabla\ge0$.
Then
\begin{equation}\label{eq:R2decrease}
 \mathcal L(X+\eta\widehat\Phi)
 \le\mathcal L(X)-\eta v^\star+\eta\nf{G}\,e_{\rm dir}
                           +\tfrac12L_\nabla q\eta^2.
\end{equation}
For rank-one normals at a nonzero boundary, the exact pseudo-Huber
direction satisfies $\nf{\Phi_\delta-\bar\Phi}\asymp\delta^{2/3}$.
Algorithm~\ref{alg:recovery} with $\epsin=O(\delta)$ gives
$e_{\rm dir}=O(\delta)$.
\end{proposition}
\begin{proof}
The reference $\bar\Phi$ attains $v^\star$.
Feasibility gives the nonnegative gap, and
$\ip{G}{\bar\Phi-\widehat\Phi}\le\nf{G}\,e_{\rm dir}$ by Cauchy--Schwarz.
Subtracting the updates gives the second formula.
The gradient Lipschitz bound gives the descent inequality
\[
 \mathcal L(X+\eta\widehat\Phi)
 \le\mathcal L(X)-\eta\ip{G}{\widehat\Phi}
                   +\tfrac12L_\nabla\eta^2\nf{\widehat\Phi}^2.
\]
Insert \eqref{eq:R2applicationbound} and
$\nf{\widehat\Phi}^2\le q\no{\widehat\Phi}^2\le q$.
The final rates follow from
Corollary~\ref{cor:rankone} and Proposition~\ref{prop:R2recovery}.
\end{proof}
The value-gap bound uses only the error norm through Cauchy--Schwarz.
For a particular computed direction, the actual gap is the pairing
$\ip{G}{\bar\Phi-\widehat\Phi}$; it can be much smaller when the
error has a small component along $G$. Section~\ref{sec:numerics}
therefore reports the gap separately from the direction error.

Under the same nonzero boundary assumptions, the center has an
additional interpretation for the quadratic objective
$\mathcal L(X)=\tfrac12\nf{X-Y}^2$, with fixed data
$Y\in\R^{m\times n}$ and $G=Y-X$.
Every exact optimal direction $\Phi\in\cF$ satisfies
$\ip{G}{\Phi}=v^\star$, so
\[
 \mathcal L(X+\eta\Phi)
 =\mathcal L(X)-\eta v^\star+\tfrac12\eta^2\nf{\Phi}^2.
\]
For $\eta>0$, the minimum-Frobenius-norm property in
Proposition~\ref{prop:cert} therefore makes $\bar\Phi$ the unique
minimizer of this one-step objective among the exact optimal directions.
This supplies a concrete reference for comparing recovery methods.
For the unrescaled SSO tangent direction introduced in
Section~\ref{sec:intro}, the same additive convention gives the
update discrepancy $\eta e_{\rm dir}$.

\paragraph{Choosing a smoothing scale}
For a prescribed total direction tolerance $\epsdir>0$, the fixed-line
bias laws give the following asymptotic scale choices.
\begin{center}
\begin{tabular}{lll}
\toprule
Method & Direction bias & Smoothing scale \\
\midrule
Pseudo-Huber, $D\ne0$ & $\asymp\delta^{1/3}$ & $\delta\asymp\epsdir^3$ \\
Pseudo-Huber, $D=0$ & $\asymp\delta^{2/3}$ & $\delta\asymp\epsdir^{3/2}$ \\
Huber or three scales & $O(\delta)$ & $\delta=O(\epsdir)$ \\
\bottomrule
\end{tabular}
\end{center}
Choose the bisection error budget $\epsin=O(\epsdir)$ as well.
The bias constant in \eqref{eq:R2finiteaccuracy} determines the
instance-dependent scale factor; Algorithm~\ref{alg:oracle} enforces
the bisection budget through its explicit width threshold.

\section{Structured numerical verification}\label{sec:numerics}
The tests check the boundary error classification, recovery of the same
center, and the effect of stopping bisection at finite accuracy.
The structured spectral calculations use 90 decimal digits to resolve
the asymptotic powers and coefficients. The algorithm tests separately
check the prescribed stopping rule and the resulting feasible updates.
Appendix~\ref{app:testinputs} specifies the ten inputs, and
Appendix~\ref{app:closed_form_check} gives an independent 120-digit check.

\subsection{Classification, coefficients, and rank protection}
For $e_{\ph}(\delta)=\nf{\Phi_\delta-\bar\Phi}$, the measured exponent
between two positive scales $\delta_2<\delta_1$ is
\begin{equation}\label{eq:observedslope}
 \frac{\log(e_{\ph}(\delta_2)/e_{\ph}(\delta_1))}
      {\log(\delta_2/\delta_1)}.
\end{equation}
Cases P01--P10 include nonsymmetric, tall and wide matrices, repeated
singular values, and empty residual blocks.
The five cases with $D\ne0$ approach exponent $1/3$; the five with
$D=0$ approach $2/3$, including P03 with its cubic residual mode.
The running example also recovers the displacement constants
$c_\lambda=4^{-1/3}$ and $\nu_\lambda=1/12$ and the matrix coefficient
$c_\lambda^2\wideD$. In particular, its displacement uses the partial
curvature $\kv=2$ rather than the full contact curvature $\kf=4$.
Appendix~\ref{app:asymptotic_tests} gives the scaled coefficient checks.

Five equality constructions attain $\rank\Theta=p+2d$;
Table~\ref{tab:rankchecks} lists $p$, $d$, the normal rank and the curvature.
P09 and P10 have rank-one and rank-two normals, respectively, and
exhibit the protected $2/3$ exponent. Fifteen additional rank-one
instances satisfy \eqref{eq:rankonecurvature}.

\subsection{Recovery of the same center}
On P01, the polynomial-tail profiles with $\rho=3/2$ and $4$ approach
the predicted exponents $1/5$ and $3/5$, while Huber approaches one.
Every error uses the same center $\bar\Phi=\diag(1,-1,0)$.
For the two-scale extrapolation, the raw spectral norm exceeds one;
radial repair restores feasibility. The repaired two- and three-scale
errors follow the rates in Corollary~\ref{cor:extrap}.
Appendix~\ref{app:asymptotic_tests} gives their values and the secondary
saturation-edge check.

\subsection{Bisection accuracy and a constrained update}\label{sec:numapp}
Use $\mathcal L(X)=\tfrac12\nf{X-Y}^2$ at $X_0=0$ with
$\ip{\Theta}{X}=0$, so the negative gradient is $G=Y$.
The rank-three test takes the $3\times3$ data in \eqref{eq:slowexample}.
The rank-one test takes $Y=\diag(1,0)$ and
$\Theta=\left[\begin{smallmatrix}1&1\\1&1\end{smallmatrix}\right]$.
Both have $G=A^\star$, $\lambda^\star=0$ and $v^\star=1$;
their centers are $\diag(1,-1,0)$ and $\diag(1,-1)$.
Algorithm~\ref{alg:oracle} starts with a doubling bracket and uses
$\epsin=\delta^3$ as its bisection direction-error budget.
All final directions are tangent-projected and radially repaired by
\eqref{eq:R2repair}.

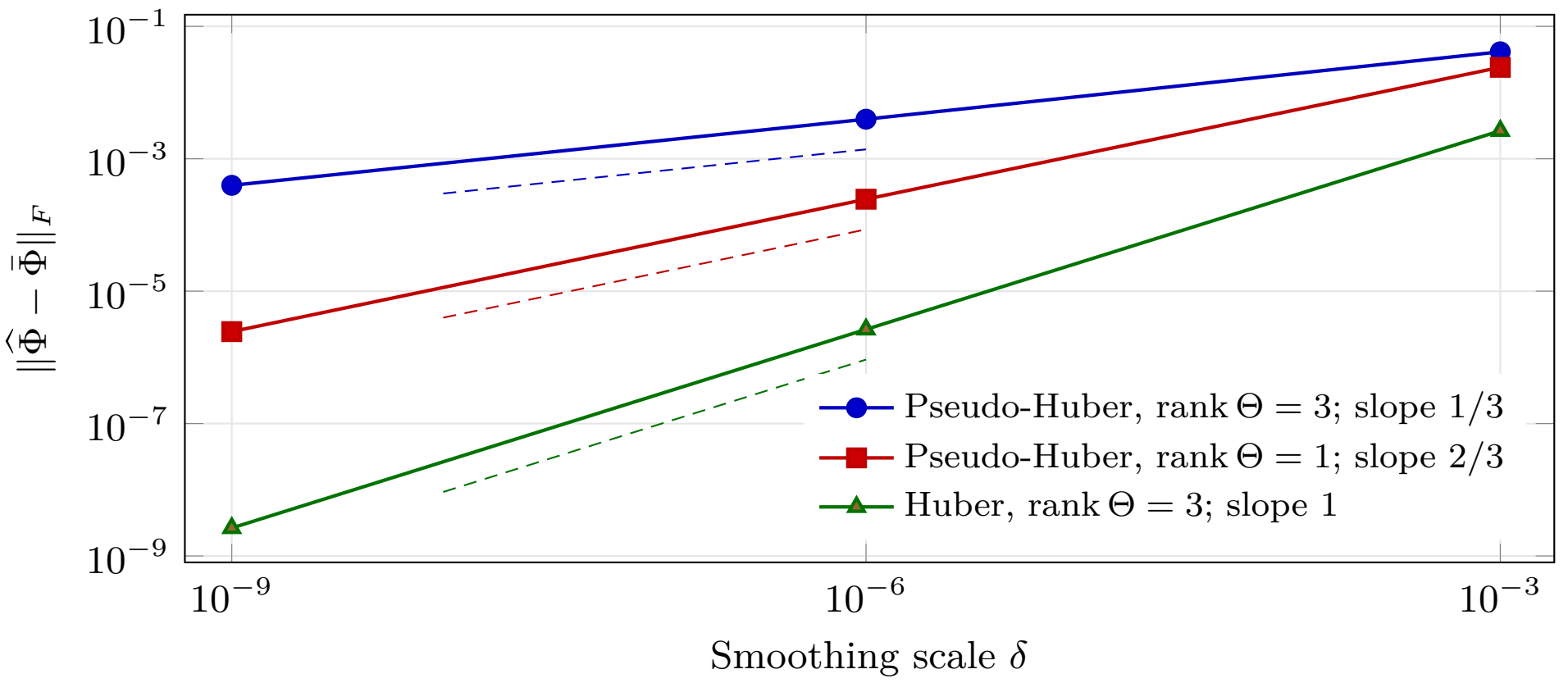
\begin{figure}[htbp]
\centering
\begin{tikzpicture}
\begin{loglogaxis}[
width=0.97\linewidth,height=6.0cm,
xlabel={Smoothing scale $\delta$},ylabel={$\|\widehat\Phi-\bar\Phi\|_F$},
xmin=6e-10,xmax=1.8e-3,ymin=8e-10,ymax=0.15,
xtick={1e-9,1e-6,1e-3},ytick={1e-9,1e-7,1e-5,1e-3,1e-1},
log basis x=10,log basis y=10,grid=major,
grid style={gray!20},tick label style={font=\small},
label style={font=\small},
legend style={at={(0.98,0.04)},anchor=south east,
  draw=none,fill=white,font=\footnotesize,
  cells={anchor=west},row sep=0pt},
mark size=2pt]
\addplot+[color=blue!75!black,mark=*,thick] coordinates {(1.000000000000e-09,3.967478029153783e-04) (1.000000000000e-06,3.960554716481276e-03) (1.000000000000e-03,4.128037889983200e-02)};
\addlegendentry{Pseudo-Huber, $\rank\Theta=3$; slope $1/3$}
\addplot[color=blue!75!black,densely dashed,no marks,forget plot] coordinates {(1.000000000000e-08,2.986464765534828e-04) (1.000000000000e-06,1.386194150768446e-03)};
\addplot+[color=red!75!black,mark=square*,thick] coordinates {(1.000000000000e-09,2.449485252060624e-06) (1.000000000000e-06,2.449040756059714e-04) (1.000000000000e-03,2.405430340156387e-02)};
\addlegendentry{Pseudo-Huber, $\rank\Theta=1$; slope $2/3$}
\addplot[color=red!75!black,densely dashed,no marks,forget plot] coordinates {(1.000000000000e-08,3.978604079236279e-06) (1.000000000000e-06,8.571642646208998e-05)};
\addplot+[color=green!45!black,mark=triangle*,thick] coordinates {(1.000000000000e-09,2.645751310308662e-09) (1.000000000000e-06,2.645750555114144e-06) (1.000000000000e-03,2.644973933699291e-03)};
\addlegendentry{Huber, $\rank\Theta=3$; slope $1$}
\addplot[color=green!45!black,densely dashed,no marks,forget plot] coordinates {(1.000000000000e-08,9.260126942899501e-09) (1.000000000000e-06,9.260126942899502e-07)};
\end{loglogaxis}
\end{tikzpicture}
\caption{Final direction errors for the constrained-update tests.
Markers show the 90-digit computations at three smoothing
scales; dashed guides have the indicated slopes.
Pseudo-Huber shows the two boundary exponents, and Huber gives linear
recovery on the same rank-three line.}
\label{fig:boundary_recovery}
\end{figure}

Figure~\ref{fig:boundary_recovery} compares the rates, and
Table~\ref{tab:R2application} includes three-scale recovery.
All 18 runs pass the feasibility checks in Appendix~\ref{app:update_checks}
and satisfy the quadratic identity and descent bound \eqref{eq:R2decrease}
at $\eta=0.1$. Coarse solves with $\epsin=\delta/100$ satisfy
\eqref{eq:R2repairLip}. The rank-three Huber value gap at $\delta=10^{-6}$
is about $10^{-12}$, well below the bound \eqref{eq:R2applicationbound}.

\section{Concluding remarks}\label{sec:discussion}
At a boundary nuclear-norm minimum, the positive and first-order
singular blocks determine the multiplier displacement, while a separate
quadratic block can dominate the direction error. Their separation gives
the two sharp pseudo-Huber rates and the rank condition that excludes
the slower one. Huber smoothing and extrapolation improve approximation
to the same minimum-norm direction, and the bisection tolerance controls
the additional error of its computation.

Uniform bounds as the boundary condition or spectral gaps vary would
quantify robustness to perturbed data. Approximate spectral evaluations
also raise the question of how to balance their error against smoothing
bias and bisection error.

\label{page:end-main}
\clearpage
\appendix
\titleformat{\section}[runin]
  {\normalfont\normalsize\bfseries}{Appendix~\thesection.}{0.5em}{}[.]
\titlespacing*{\section}{\parindent}{9pt plus 2pt minus 1pt}{0.5em}
\section{Construction of the orthogonal normal form}\label{app:normalform}

The purpose of the construction is to separate singular values of orders
$1$, $|s|$, and $s^2$ by transformations that preserve singular values.
We use the fixed bases $U_\star=[U_+,U_0]$ and
$V_\star=[V_+,V_0]$ from Section~\ref{sec:geometry}; matrices without a
moving argument denote these fixed bases. All matrix remainders below
are measured in the Frobenius norm.

\begin{proof}[\hypertarget{proof:normalform}{}Proof of Lemma~\ref{lem:normalform}]
\emph{Separating the positive singular values of $A^\star$.}
Define the left and right Gram matrices
\[
 M_L(s)=A(s)A(s)^T\in\R^{m\times m},
 \qquad M_R(s)=A(s)^TA(s)\in\R^{n\times n}.
\]
At $s=0$, their $r$ positive eigenvalues are
$\sigma_1^2,\ldots,\sigma_r^2$, separated from zero by
$\sigma_r^2>0$. Choose a positively oriented closed contour
$\Gamma$ in the complex plane that encloses this positive cluster and
excludes zero. The same contour separates the clusters for sufficiently
small real $s$. With $\mathrm i^2=-1$ and complex integration variable $z$,
the spectral projectors are
\[
 \Pi_L(s)=\frac{1}{2\pi\mathrm i}
       \int_\Gamma(zI_m-M_L(s))^{-1}\,dz,\qquad
 \Pi_R(s)=\frac{1}{2\pi\mathrm i}
       \int_\Gamma(zI_n-M_R(s))^{-1}\,dz.
\]
For real $s$, they are real orthogonal projectors of rank $r$, analytic
in $s$~\cite{kato1995}. Thus repeated eigenvalues within the selected
cluster cause no loss of analyticity. The resolvent identity
\[
 (zI_m-M_L(s))^{-1}A(s)
       =A(s)(zI_n-M_R(s))^{-1}
\]
implies $\Pi_L(s)A(s)=A(s)\Pi_R(s)$.

We next turn these projectors into orthonormal bases. For any analytic
orthogonal projector $\Pi(s)\in\R^{N\times N}$ of constant rank $k$,
let $Q_{\rm ref}\in\R^{N\times k}$ be an orthonormal basis of
$\operatorname{range}\Pi(0)$. Then
\begin{equation}\label{eq:R2frame}
 Q(s)=\Pi(s)Q_{\rm ref}
       [Q_{\rm ref}^T\Pi(s)Q_{\rm ref}]^{-1/2}
\end{equation}
is an analytic orthonormal basis of $\operatorname{range}\Pi(s)$:
the bracketed Gram matrix equals $I_k$ at zero and remains positive
definite nearby, and $\Pi(s)^2=\Pi(s)$ gives $Q(s)^TQ(s)=I_k$.
Apply this formula to $\Pi_L$, $I_m-\Pi_L$, $\Pi_R$, and $I_n-\Pi_R$,
with reference bases $U_+$, $U_0$, $V_+$, and $V_0$, respectively.
The resulting matrices have sizes
\[
 \begin{array}{c|cccc}
 \text{basis}&U_+(s)&U_0(s)&V_+(s)&V_0(s)\\ \hline
 \text{size}&m\times r&m\times m_0&n\times r&n\times n_0 .
 \end{array}
\]
They equal their fixed counterparts at zero. The intertwining identity
makes the two off-diagonal blocks vanish, so
\[
 [U_+(s),U_0(s)]^T A(s)[V_+(s),V_0(s)]
       =\diag(A_+(s),B_0(s)),
 \quad
 B_0(s)=U_0(s)^TA(s)V_0(s).
\]
Here $A_+(s)\in\R^{r\times r}$, $A_+(0)=\Sigma$, and
$B_0(s)\in\R^{m_0\times n_0}$, $B_0(0)=0$.

\emph{Computing the first low-spectrum block.}
Differentiating the projector identity gives
\[
 \Pi^2=\Pi,\qquad \Pi(0)\Pi'(0)\Pi(0)=0.
\]
Consequently the basis choice \eqref{eq:R2frame} has zero derivative
within its own subspace at zero. Differentiating
$U_0(s)^TA(s)V_+(s)=0$ and
$U_+(s)^TA(s)V_0(s)=0$ therefore determines the remaining components:
\[
 U_0'(0)^TU_+\Sigma+\Theta_{0+}=0,\qquad
 \Sigma V_+^TV_0'(0)+\Theta_{+0}=0.
\]
In the original matrix coordinates, this gives
\begin{equation}\label{eq:framederiv}
 U_0'(0)=-U_+\Sigma^{-1}\Theta_{0+}^T,\qquad
 V_0'(0)=-V_+\Sigma^{-1}\Theta_{+0} .
\end{equation}
The coefficient of $s$ in $B_0(s)$ is $C$. Its coefficient of $s^2$
is the sum
\[
 \begin{aligned}
 &U_0'(0)^TA^\star V_0'(0)
  +U_0'(0)^T\Theta V_0
  +U_0^T\Theta V_0'(0)\\
 &\hspace{2em}
   =\Theta_{0+}\Sigma^{-1}\Theta_{+0}
    -\Theta_{0+}\Sigma^{-1}\Theta_{+0}
    -\Theta_{0+}\Sigma^{-1}\Theta_{+0}\\
 &\hspace{2em}=-\Theta_{0+}\Sigma^{-1}\Theta_{+0}.
 \end{aligned}
\]
The two terms involving second derivatives of the bases vanish because
$A^\star V_0=0$ and $U_0^TA^\star=0$. Hence
\begin{equation}\label{eq:lowerjet}
 B_0(s)=sC-s^2\Theta_{0+}\Sigma^{-1}\Theta_{+0}+O(s^3).
\end{equation}
This computation uses orthogonal bases throughout, so it identifies a
coefficient of the actual low-spectrum block.

\emph{Separating the first-order and residual singular values.}
The quotient $\widetilde B(s)=B_0(s)/s\in\R^{m_0\times n_0}$ has an analytic extension
at zero, with value $C$. Its $p$ positive squared singular values
are separated from zero by $\sigma_{C,\min}^2$. Apply the preceding projector
and basis construction to $\widetilde B(s)$, using the fixed reference
bases $[U_C,U_{C,\perp}]$ and $[V_C,V_{C,\perp}]$. It gives analytic orthogonal
matrices
\[
 \mathcal U_C(s)\in\R^{m_0\times m_0},\quad
 \mathcal V_C(s)\in\R^{n_0\times n_0},\qquad
 \mathcal U_C(0)=[U_C,U_{C,\perp}],\quad \mathcal V_C(0)=[V_C,V_{C,\perp}],
\]
such that
\[
 \mathcal U_C(s)^T\widetilde B(s)\mathcal V_C(s)
       =\diag(A_{\mathrm{act}}(s),B_{\rm res}(s)).
\]
The active block $A_{\mathrm{act}}(s)$ is $p\times p$, with $A_{\mathrm{act}}(0)=\Sigma_C$, and
$B_{\rm res}(s)$ is $(m_0-p)\times(n_0-p)$, with $B_{\rm res}(0)=0$.
In its derivative at zero, every basis-derivative term has $C$ on
a residual side and vanishes. Thus
\[
 B_{\rm res}'(0)
       =-U_{C,\perp}^T\Theta_{0+}\Sigma^{-1}\Theta_{+0}V_{C,\perp}=D_0.
\]
Writing $B_{\rm res}(s)=sA_{\mathrm{res}}(s)$ gives analytic
$A_{\mathrm{res}}$ with $A_{\mathrm{res}}(0)=D_0$.
Thus the second reduction of the actual low-spectrum block is
\[
 \mathcal U_C(s)^TB_0(s)\mathcal V_C(s)
   =\diag(sA_{\mathrm{act}}(s),s^2A_{\mathrm{res}}(s)).
\]
The complete transformations in the statement are
\[
 \begin{aligned}
 \cU(s)&=[U_+(s),U_0(s)]\diag(I_r,\mathcal U_C(s)),\\
 \cV(s)&=[V_+(s),V_0(s)]\diag(I_r,\mathcal V_C(s)).
 \end{aligned}
\]
They give \eqref{eq:normalform}, with initial bases
$[U_+,U_0U_C,U_0U_{C,\perp}]$ and $[V_+,V_0V_C,V_0V_{C,\perp}]$.
If $D_0=0$, analyticity and $A_{\mathrm{res}}(0)=0$ give
$A_{\mathrm{res}}(s)=O(s)$, so the actual residual block is $O(s^3)$.
The construction also applies when either residual dimension is zero;
the corresponding empty blocks contribute zero.

\emph{Computing the trace coefficients.}
For any $A\in\R^{m\times n}$, introduce the real symmetric dilation
\[
 \mathcal H(A)=
 \begin{pmatrix}0_{m\times m}&A\\A^T&0_{n\times n}\end{pmatrix}
       \in\R^{(m+n)\times(m+n)}.
\]
Its spectrum consists of the positive and negative singular values of
$A$ and any additional zeros. Along our line,
$d\mathcal H(A(s))/ds=\mathcal H(\Theta)$.
Let $(\xi_\ell,w_\ell)$, $1\le\ell\le m+n$, be an orthonormal
eigensystem of $\mathcal H(A^\star)$, with
$w_\ell\in\R^{m+n}$ and real dilation eigenvalues $\xi_\ell$.
Let $\mathcal I_+$ be the $r$ indices of its positive eigenvalues.
Their sum along the separated cluster is $\ell_+(s)=\nn{A_+(s)}$.

For completeness, if $\Pi_{\rm d}(s)$ is this cluster projector,
its derivative has, between the selected cluster and its complement,
the entries
\[
 w_j^T\Pi_{\rm d}'(0)w_i
   =\frac{w_j^T\mathcal H(\Theta)w_i}
          {\xi_i-\xi_j},
 \qquad i\in\mathcal I_+,\quad j\notin\mathcal I_+ .
\]
This follows by differentiating
\[
 \mathcal H(A(s))\Pi_{\rm d}(s)
   =\Pi_{\rm d}(s)\mathcal H(A(s)).
\]
The dilation is affine in $s$, so its second derivative vanishes.
Differentiation of the cluster trace yields
\begin{equation}\label{eq:R2clusterhessian}
 \ell_+''(0)=
  2\sum_{i\in\mathcal I_+}\sum_{j\notin\mathcal I_+}
  \frac{|w_j^T\mathcal H(\Theta)w_i|^2}
       {\xi_i-\xi_j}.
\end{equation}
The factor two comes from the two off-diagonal blocks of
$\Pi_{\rm d}'(0)$. Every denominator separates the selected cluster from
its complement; none separates eigenvalues within a repeated cluster.

Let $u_i$ and $v_i$, $1\le i\le r$, denote the columns of the fixed
$U_+$ and $V_+$, so they are singular vectors of $A^\star$.
The dilation eigenvectors associated with $\pm\sigma_i$ are
$(u_i,\pm v_i)/\sqrt2$. The matrix element from a positive vector $i$
to a negative vector $j$ is $((\Theta_{++})_{ji}-(\Theta_{++})_{ij})/2$.
Their contribution to \eqref{eq:R2clusterhessian} is therefore
\[
 \sum_{i,j=1}^r\frac{((\Theta_{++})_{ji}-(\Theta_{++})_{ij})^2}{2(\sigma_i+\sigma_j)}
       =\sum_{1\le i<j\le r}
             \frac{((\Theta_{++})_{ij}-(\Theta_{++})_{ji})^2}{\sigma_i+\sigma_j}.
\]
For the zero eigenvalues, use the vectors
$(U_0{}_{:,\alpha},0)$, $1\le\alpha\le m_0$, and
$(0,V_0{}_{:,\beta_0})$, $1\le\beta_0\le n_0$.
Their matrix elements are $(\theta_{0+,i})_\alpha/\sqrt2$ and
$(\theta_{+0,i})_{\beta_0}/\sqrt2$, where $\theta_{0+,i}=(\Theta_{0+})_{:,i}$ and
$\theta_{+0,i}=(\Theta_{+0})_{i,:}^T$. Thus their total contribution is
$\sum_{i=1}^r(\|\theta_{0+,i}\|^2+\|\theta_{+0,i}\|^2)/\sigma_i$.
We have proved that the scalar positive-cluster curvature is
\begin{equation}\label{eq:kplus}
 \kappa_+:=\ell_+''(0)=
 \sum_{1\le i<j\le r}\frac{((\Theta_{++})_{ij}-(\Theta_{++})_{ji})^2}{\sigma_i+\sigma_j}
 +\sum_{i=1}^r\frac{\|\theta_{0+,i}\|^2+\|\theta_{+0,i}\|^2}{\sigma_i}.
\end{equation}
Also $\ell_+'(0)=\tr \Theta_{++}=a$, so
\begin{equation}\label{eq:lplus}
 \ell_+(s)=\nn{A^\star}+as+\tfrac12\kappa_+s^2+O(s^3).
\end{equation}

The active sum $\ell_{\mathrm{act}}(s)=\nn{A_{\mathrm{act}}(s)}$ is the sum of the $p$ positive
singular values of $\widetilde B(s)$ separated from zero. Its first
derivative is the pairing with the corresponding partial polar factor
$\msign(C)=U_CV_C^T$~\cite{lewis2001}. Applying this derivative to
$\widetilde B'(0)=-\Theta_{0+}\Sigma^{-1}\Theta_{+0}$ gives
\begin{equation}\label{eq:gjet}
 \ell_{\mathrm{act}}(s)=\beta+s\ip{\msign(C)}{-\Theta_{0+}\Sigma^{-1}\Theta_{+0}}+O(s^2).
\end{equation}
In particular,
$\kv=\kappa_+-2\ip{\msign(C)}{-\Theta_{0+}\Sigma^{-1}\Theta_{+0}}$, as required for the scalar
balance in the main text.

Finally, an alternative analytic choice of bases within the first
low-spectrum spaces, with the same bases at zero, changes $\widetilde B'(0)$ by
$\Xi_L^TC+C\Xi_R$, where
$\Xi_L\in\R^{m_0\times m_0}$ and
$\Xi_R\in\R^{n_0\times n_0}$ are skew-symmetric.
Its residual compression vanishes because
$U_{C,\perp}^TC=0$ and $CV_{C,\perp}=0$.
Its pairing with $\msign(C)$ is zero as well:
\[
 \ip{\msign(C)}{\Xi_L^TC+C\Xi_R}
   =\tr(\Sigma_CU_C^T\Xi_L^TU_C)+\tr(\Sigma_CV_C^T\Xi_RV_C)=0.
\]
Each trace pairs the symmetric diagonal $\Sigma_C$ with a skew-symmetric
matrix. Thus $D_0$ and the coefficient in \eqref{eq:gjet} are unchanged
by this within-subspace choice.
\end{proof}

\subsection{The flat case of boundary rigidity}\label{app:flat_rigidity}
We complete the zero-curvature step in Theorem~\ref{thm:rigidity}.
If $\kv=0$, every square in \eqref{eq:kv} and \eqref{eq:squares}
vanishes. Hence $\Theta_{++}=\Theta_{++}^T$,
$P_{C,L}^{\perp}\Theta_{0+}=0$, $\Theta_{+0}P_{C,R}^{\perp}=0$, and
$U_C^T\Theta_{0+}=-V_C^T\Theta_{+0}^T$.
Put $\Theta_{+a}=\Theta_{+0}V_C\in\R^{r\times p}$, the coupling to
the active right singular space of $C$.
The only nonzero part of the line, between the positive and active
coordinates, is the $(r+p)\times(r+p)$ matrix
\[
 A_c(s)=\begin{pmatrix}\Sigma+s\Theta_{++}&s\Theta_{+a}\\-s\Theta_{+a}^T&s\Sigma_C\end{pmatrix}.
\]
Left multiplication by $\diag(I_r,-I_p)$ preserves its singular values
and gives the symmetric matrix
\[
 \begin{pmatrix}
 \Sigma+s\Theta_{++}&s\Theta_{+a}\\
 s\Theta_{+a}^T&-s\Sigma_C
 \end{pmatrix}.
\]
For small negative $s$, $\Sigma+s\Theta_{++}\succeq(\sigma_r/2)I_r$.
Its Schur complement satisfies
\[
 \begin{aligned}
 &-s\Sigma_C-s^2\Theta_{+a}^T
                  (\Sigma+s\Theta_{++})^{-1}\Theta_{+a}\\
 &\qquad\succeq |s|\left(\sigma_{C,\min}
              -\frac{2|s|\no{\Theta_{+a}}^2}{\sigma_r}\right)I_p
              \succ0.
 \end{aligned}
\]
Thus the displayed symmetric matrix is positive definite.
Its nuclear norm equals its trace,
and $\phi(s)=\tr\Sigma+s(a-\beta)=\phi(0)$ on a negative interval.

\subsection{Sharp rank constructions and the rank-one bound}\label{app:sharp_rank}
\paragraph{Sharpness in Theorem~\ref{thm:rank}}
For equality, take $r=d$, $m_0=n_0=p+d$ and
\begin{equation}\label{eq:sharpconstruction}
 \begin{aligned}
 \Sigma&=I_d,\qquad C=\diag(I_p,0_d),
                     \qquad\Theta_{++}=(p/d)I_d,\\
 \Theta_{0+}&=\begin{pmatrix}0_{p\times d}\\I_d\end{pmatrix},
 \qquad\Theta_{+0}=\begin{pmatrix}0_{d\times p}&I_d\end{pmatrix}.
 \end{aligned}
\end{equation}
Then $a=\beta=p$, $D_0=-I_d$, and $\kv=2d$.
After the $\Sigma_C$ pivot the remaining block is
$\left[\begin{smallmatrix}(p/d)I_d&I_d\\I_d&0\end{smallmatrix}\right]$.
Its kernel equations force both components to be zero, so it has rank
$2d$. The minimizer is unique by Theorem~\ref{thm:rigidity}.

\begin{proof}[Proof of Corollary~\ref{cor:rankone}]
Factor $\Theta=uv^T$ with nonzero $u\in\R^m$, $v\in\R^n$.
These factors need not be unit vectors.
Their fixed singular-coordinate components are
\[
 \begin{aligned}
 u_+&=U_+^Tu\in\R^r,&u_0&=U_0^Tu\in\R^{m_0},\\
 v_+&=V_+^Tv\in\R^r,&v_0&=V_0^Tv\in\R^{n_0}.
 \end{aligned}
\]
The blocks are $\Theta_{++}=u_+v_+^T$, $\Theta_{+0}=u_+v_0^T$,
$\Theta_{0+}=u_0v_+^T$, and $C=u_0v_0^T$.
The boundary gives $u_+^Tv_+=\beta=\|u_0\|\|v_0\|>0$, hence $p=1$,
$U_C=u_0/\|u_0\|$, and $V_C=v_0/\|v_0\|$.
Since $P_{C,L}^{\perp}\Theta_{0+}=0$ and
$\Theta_{+0}P_{C,R}^{\perp}=0$, the residual block is zero.

The $r$ active scalar squares in \eqref{eq:squares} sum to
\[
 \begin{aligned}
 \sum_{i=1}^r|U_C^T\theta_{0+,i}+V_C^T\theta_{+0,i}|^2
   &=\bigl\|\|u_0\|v_++\|v_0\|u_+\bigr\|^2\\
   &=\bigl\|\|u_0\|v_+-\|v_0\|u_+\bigr\|^2
                          +4\beta\,u_+^Tv_+\\
   &\ge4\beta^2 .
 \end{aligned}
\]
Their denominators satisfy $\sigma_i\le\sigma_1$, so retaining
these squares in \eqref{eq:kv} proves
$\kv\ge4\beta^2/\sigma_1>0$.
Since the only positive singular value of $C$ is $\sigma_{C,1}=\beta$,
$b_{\ph}=1/(2\beta)$ and
$c_{\lambda}^3=1/(2\beta\kv)\le\sigma_1/(8\beta^3)$.
This gives the displacement coefficient bound.
For any $\sigma>0$ and $\beta>0$, take
\[
 A^\star=\diag(\sigma,0),\qquad
 \Theta=\beta\begin{pmatrix}1&1\\1&1\end{pmatrix}.
\]
Then $\kv=4\beta^2/\sigma$ and both inequalities are equalities.
Theorems~\ref{thm:rigidity} and~\ref{thm:PH} supply uniqueness and
the exact direction rate.
\end{proof}

\section{Proof of the pseudo-Huber classification}\label{app:ph}

We prove Theorem~\ref{thm:PH} in two stages. First we locate the
minimizing multiplier by comparing the three spectral blocks.
We then reparameterize the direction by $t=\delta^{1/3}$, identify its
analytic branch, and compute the first two matrix coefficients.
Throughout, $\kv>0$ and $a=\beta>0$; all constants depend on the fixed
matrices and all matrix estimates use the Frobenius norm.

\subsection{Scalar balance and root localization}\label{sec:root}

\begin{proof}[Proof of Theorem~\ref{thm:PH}, scalar part]
The scalar function to be solved is the tangency residual
\[
 h_\delta(\lambda^\star+s)
   =\ip{\Theta}{Z_\delta(A(s))}
   =\frac{d}{ds}\Psi_\delta(A(s)).
\]
Its unique zero is $s_\delta$, by Proposition~\ref{prop:PHselection}.
Orthogonal invariance of $\Psi_\delta$ gives, from
\eqref{eq:normalform},
\[
 \Psi_\delta(A(s))
   =\Psi_\delta(A_+(s))+
     \Psi_\delta(sA_{\mathrm{act}}(s))+\Psi_\delta(s^2A_{\mathrm{res}}(s)).
\]
Each spectral function on the right uses the dimensions of its own
matrix input. We differentiate this scalar identity, so moving bases
introduce no additional terms.

\emph{Positive and first-order singular values.}
The $r\times r$ block $A_+(s)$ remains uniformly invertible near zero.
On such matrices,
$Z_\delta(A_+(s))=\msign(A_+(s))+O(\delta^2)$.
Since $A_+'(s)$ is bounded, the first block contributes
$\ell_+'(s)+O(\delta^2)$ to the scalar derivative.

For the active block $A_{\mathrm{act}}(s)\in\R^{p\times p}$, let $k_j(s)$,
$1\le j\le p$, be positive analytic singular-value branches, labeled so
that $k_j(0)=\sigma_{C,j}$. These branches exist by applying the analytic
eigenvalue theorem to the positive definite analytic matrix
$A_{\mathrm{act}}(s)^TA_{\mathrm{act}}(s)$ and taking positive square roots
\cite{kato1995,bunsegerstner1991}.
They are local analytic branches, rather than singular values required
to remain in decreasing order. They satisfy
$k_j(s)\ge \sigma_{C,\min}/2$ and have bounded derivatives on a fixed small
interval. For $s<0$, the active singular values of $A(s)$ are therefore
\[
 \sigma_{a,j}(s)=-s\,k_j(s)>0,\qquad
 \sigma_{a,j}'(s)=-k_j(s)-s\,k_j'(s).
\]
The subscript $a$ identifies this active cluster; these quantities are
singular values, not scalar multipliers.

For a positive scalar $u$ in a region where $\delta/u$ is small,
\[
 \frac{u}{\sqrt{u^2+\delta^2}}
 =1-\frac{\delta^2}{2u^2}
      +O(\delta^4/u^4).
\]
The estimate is uniform, for example when $0\le\delta/u\le1/2$.
The unsmoothed active derivative is
$\sum_{j=1}^p\sigma_{a,j}'(s)=-\ell_{\mathrm{act}}(s)-s\ell_{\mathrm{act}}'(s)$.
Subtracting it from the smoothed derivative gives
\begin{equation}\label{eq:activeestimate}
 \begin{aligned}
 &\sum_{j=1}^p\sigma_{a,j}'(s)
  \left(\frac{\sigma_{a,j}(s)}
   {\sqrt{\sigma_{a,j}(s)^2+\delta^2}}-1\right)\\
 &\qquad=\frac{\delta^2}{2s^2}\sum_{j=1}^p \sigma_{C,j}^{-1}
       +O\left(\frac{\delta^2}{|s|}
                    +\frac{\delta^4}{s^4}\right).
 \end{aligned}
\end{equation}
Indeed, $-\sigma_{a,j}'(s)=\sigma_{C,j}+O(|s|)$ and
$\sigma_{a,j}(s)^2=s^2(\sigma_{C,j}^2+O(|s|))$ give the first remainder;
multiplying the fourth-order scalar remainder by a bounded derivative
gives the second.

\emph{Residual singular values.}
Write the residual block as
$R_0(s)=s^2A_{\mathrm{res}}(s)\in
 \R^{(m_0-p)\times(n_0-p)}$.
It satisfies
\[
 R_0(s)=s^2D_0+O(s^3),\qquad
 R_0'(s)=2sD_0+O(s^2).
\]
For any matrix $X$ of these dimensions with $\no{X}/\delta$ small,
the inverse-square-root expansion gives
\[
 Z_\delta(X)=\frac{X}{\delta}
          +O\left(\frac{\nf{X}^3}{\delta^3}\right).
\]
This is a matrix estimate and permits arbitrary rank of $X$.
The residual contribution is consequently
\begin{equation}\label{eq:residualestimate}
 \ip{R_0'(s)}{Z_\delta(R_0(s))}
 =\frac{2s^3}{\delta}\nf{D}^2
   +O\left(\frac{s^4}{\delta}
                     +\frac{|s|^7}{\delta^3}\right).
\end{equation}
Here $\nf{D_0}=\nf{D}$ by the fixed orthonormal compression.
The first remainder follows from
$\ip{R_0'}{R_0}=2s^3\nf{D_0}^2+O(s^4)$;
the second is the product of $O(|s|)$ with the cubic matrix remainder.

\emph{Locating the true zero.}
For now $x>0$ is a test displacement, with $s=-x$.
Fix $0<z_-<z_+<\infty$ and put
$x=z\delta^{2/3}$, where $z\in[z_-,z_+]$.
Then $\delta/x\to0$ and $x^2/\delta\to0$ uniformly in $z$,
so both preceding expansions apply. Equations
\eqref{eq:lplus}--\eqref{eq:gjet} imply
\[
 \ell_+'(s)-\ell_{\mathrm{act}}(s)-s\ell_{\mathrm{act}}'(s)=\kv s+O(s^2),
\]
using $a=\beta$ and
$\kv=\kappa_+-2\ip{\msign(C)}{-\Theta_{0+}\Sigma^{-1}\Theta_{+0}}$.
Combining the three block derivatives yields
\begin{equation}\label{eq:balance}
 h_\delta(\lambda^\star-x)
 =-\kv x+b_{\ph}\frac{\delta^2}{x^2}
      -2\nf{D}^2\frac{x^3}{\delta}+O(\delta^{4/3}),
\end{equation}
uniformly on this test interval. The two active remainders are
$O(\delta^{4/3})$, the two residual remainders are
$O(\delta^{5/3})$, and the remaining Taylor remainder is
$O(x^2)=O(\delta^{4/3})$.

After division by $\delta^{2/3}$, the residual converges uniformly
to the scalar expression
\[
 -\kv z+b_{\ph}z^{-2}\qquad(z>0).
\]
This is strictly decreasing and has the unique zero
$c_\lambda=(b_{\ph}/\kv)^{1/3}$. Choose arbitrary fixed
$0<z_-<c_\lambda<z_+$. For small $\delta$,
\[
 h_\delta(\lambda^\star-z_+\delta^{2/3})<0
       <h_\delta(\lambda^\star-z_-\delta^{2/3}).
\]
Monotonicity of $h_\delta$ puts its actual zero between these test
points, and hence
\[
 -z_+\delta^{2/3}\le s_\delta\le-z_-\delta^{2/3}.
\]
Taking $z_-\uparrow c_\lambda$ and $z_+\downarrow c_\lambda$ proves
$-s_\delta/\delta^{2/3}\to c_\lambda$. In particular the minimizing
displacement is negative for all sufficiently small $\delta$.
\end{proof}

\subsection{Analytic continuation and matrix coefficients}\label{sec:analytic}

\begin{proof}[Proof of Theorem~\ref{thm:PH}, matrix and analytic parts]
Set $t=\delta^{1/3}$ and $s=-t^2z$, so the positive displacement
is $x=-s=t^2z$. Initially $t>0$ and $z$ is near
the positive constant $c_\lambda$. Define a two-variable matrix field by
$\mathcal Z(t,z)=Z_{t^3}(A^\star-t^2z\Theta)$.
Scaling the three blocks of \eqref{eq:normalform} gives the exact formula
\begin{equation}\label{eq:analyticZ}
 \mathcal Z(t,z)=\cU(s)
 \diag\!\left(
 Z_{t^3}(A_+(s)),-Z_{t/z}(A_{\mathrm{act}}(s)),
 Z_1(tz^2A_{\mathrm{res}}(s))\right)\cV(s)^T.
\end{equation}
For example,
$Z_{t^3}(-t^2zA_{\mathrm{act}})=-Z_{t/z}(A_{\mathrm{act}})$ for $t,z>0$.
The residual block has scale $s^2=t^4z^2$, which gives its
argument $tz^2A_{\mathrm{res}}(s)$ after division by $t^3$.
Thus $\mathcal Z$ takes values in the original
$\R^{m\times n}$ matrix space.

For a uniformly invertible square matrix $X$,
\[
 Z_\delta(X)=X(X^TX+\delta^2I)^{-1/2}
\]
is analytic in $X$ and $\delta^2$ near $\delta=0$.
For a rectangular matrix $Y$ near zero,
$Z_1(Y)=Y(I+Y^TY)^{-1/2}$ is analytic in its entries.
The bases and all three blocks are analytic in $s=-t^2z$.
Consequently the right side of \eqref{eq:analyticZ} defines an analytic
extension to an open real neighborhood of $(t,z)=(0,c_\lambda)$, including
negative $t$. The original smoothed direction is obtained on $t>0$.

At $t=0$, the positive and active polar factors equal $I_r$ and
$I_p$, respectively. The initial bases of Appendix~\ref{app:normalform}
therefore give
\begin{equation}\label{eq:firstt}
 \mathcal Z(0,z)=\bar\Phi,\qquad
 \partial_t\mathcal Z(0,z)=z^2\wideD.
\end{equation}
To see the derivative, basis motion starts at order $s=O(t^2)$,
the active smoothing correction starts at $t^2$, and the positive
correction starts at $t^6$. Only the residual block has an order-$t$
term, namely $tz^2D_0$; its lift is $tz^2\wideD$.

\emph{Selecting the analytic root branch.}
The field must satisfy the scalar tangency condition.
Both its constant and linear $t$ coefficients vanish after pairing
with $\Theta$, because
\[
 \ip{\Theta}{\bar\Phi}=0,\qquad
 \ip{\Theta}{\wideD}=\ip{C}{D}=0 .
\]
The second identity uses $C^TD=0$, a consequence of $U_C^TP_{C,L}^{\perp}=0$.
These identities hold for every $z$ near $c_\lambda$, so the analytic series
of $\ip{\Theta}{\mathcal Z(t,z)}$ is divisible by $t^2$.
Define the real-valued scaled residual by
\[
 h_{\mathrm{sc}}(t,z)
   =t^{-2}h_{t^3}(\lambda^\star-t^2z)
   =t^{-2}\ip{\Theta}{\mathcal Z(t,z)}\qquad(t>0).
\]
The divisibility just proved gives its analytic continuation to $t=0$.
Equation~\eqref{eq:balance} gives
\begin{equation}\label{eq:IFTbalance}
 h_{\mathrm{sc}}(t,z)
   =-\kv z+b_{\ph}z^{-2}
      -2\nf{D}^2z^3t+O(t^2).
\end{equation}
At the point where the limiting scalar balance vanishes,
\[
 h_{\mathrm{sc}}(0,c_\lambda)=0,\qquad
 \partial_zh_{\mathrm{sc}}(0,c_\lambda)=-3\kv<0,\qquad
 \partial_th_{\mathrm{sc}}(0,c_\lambda)=-2\nf{D}^2c_\lambda^3.
\]
The analytic implicit function theorem therefore gives a unique
analytic solution $z(t)$ near $c_\lambda$, with
\[
 z(0)=c_\lambda,\qquad
 z'(0)=-\frac{\partial_th_{\mathrm{sc}}(0,c_\lambda)}
                   {\partial_zh_{\mathrm{sc}}(0,c_\lambda)}
       =-\frac{2c_\lambda^3\nf{D}^2}{3\kv}=-\nu_\lambda.
\]
The scalar localization already proved puts
$-s_\delta/\delta^{2/3}$ near $c_\lambda$, so the minimizing zero lies on
this branch. Hence
$s_\delta=-t^2z(t)=-c_\lambda t^2+\nu_\lambda t^3+O(t^4)$.
The associated matrix path is $\Phi_{t^3}=\mathcal Z(t,z(t))$.
This proves \eqref{eq:sexp} and its convergent power series.

\emph{Computing the matrix coefficients.}
Separate the motion of the large and active polar factors from their
smoothing corrections. Define the partial polar field
$\Phi_{\mathrm{part}}(s)\in\R^{m\times n}$, which retains the two
nonresidual polar blocks on the negative branch, and the fixed active
weighting matrix $M_{\mathrm{act}}\in\R^{m\times n}$ by
\begin{equation}\label{eq:W}
 \begin{aligned}
 \Phi_{\mathrm{part}}(s)&=\cU(s)
   \diag(\msign(A_+(s)),-\msign(A_{\mathrm{act}}(s)),0)\cV(s)^T,\\
 M_{\mathrm{act}}&=U_0U_C\Sigma_C^{-2}V_C^TV_0^T .
 \end{aligned}
\end{equation}
The partial field extends analytically through zero with
$\Phi_{\mathrm{part}}(0)=\bar\Phi$.
The matrix $M_{\mathrm{act}}$ uses the same fixed singular vectors as
$\wideC=U_0U_C\Sigma_CV_C^TV_0^T$, with weights $\sigma_{C,j}^{-2}$.
On the negative branch, differentiating the two nonresidual nuclear
norms gives
\[
 \ip{\Theta}{\Phi_{\mathrm{part}}(s)}
       =\ell_+'(s)-\ell_{\mathrm{act}}(s)-s\ell_{\mathrm{act}}'(s).
\]
To justify this identity in moving coordinates, differentiating the
bases adds a skew-symmetric generator on each side of the block matrix.
Its pairing with the displayed polar blocks is zero: for a square
invertible block $X$ with polar factor $P$, both $XP^T$ and $P^TX$
are symmetric, and the trace of their product with a skew-symmetric
generator vanishes. The residual polar block is zero, so cross-cluster
terms contribute zero as well. The remaining pairing with the
block derivative is exactly $\ell_+'(s)-\ell_{\mathrm{act}}(s)-s\ell_{\mathrm{act}}'(s)$.
In particular, differentiating at zero gives
$\ip{\Theta}{\Phi_{\mathrm{part}}'(0)}=\kv$.

For the active smoothing correction, evaluate at $s=-t^2z$.
Uniform invertibility of $A_{\mathrm{act}}(s)$ gives
\[
 \begin{aligned}
 -Z_{t/z}(A_{\mathrm{act}}(s))
  &=-\msign(A_{\mathrm{act}}(s))\\
  &\quad+\frac{t^2}{2z^2}\msign(A_{\mathrm{act}}(s))
      (A_{\mathrm{act}}(s)^TA_{\mathrm{act}}(s))^{-1}+O(t^4).
 \end{aligned}
\]
At $t=0$ its correction lifts to $t^2M_{\mathrm{act}}/(2z^2)$.
For the residual contribution, use
\[
 A_{\mathrm{res}}(-t^2z)=D_0+O(t^2),\qquad
 Z_1(Y)=Y+O(\nf{Y}^3).
\]
Inserting $z(t)=c_\lambda-\nu_\lambda t+O(t^2)$ in these terms gives
\[
 \begin{aligned}
 \Phi_{\mathrm{part}}(-t^2z(t))
   &=\bar\Phi-c_\lambda t^2\Phi_{\mathrm{part}}'(0)+O(t^3),\\
 \text{active smoothing correction}
   &=\frac{t^2}{2c_\lambda^2}M_{\mathrm{act}}+O(t^3),\\
 \text{residual direction}
   &=(c_\lambda^2t-2c_\lambda\nu_\lambda t^2)\wideD+O(t^3).
 \end{aligned}
\]
The bases vary by $O(t^2)$, so their multiplication changes the
leading residual term only at order $t^3$.
Consequently the coefficient in \eqref{eq:Pexp} is
\begin{equation}\label{eq:M}
 M_\Phi=-c_\lambda\Phi_{\mathrm{part}}'(0)+\frac{1}{2c_\lambda^2}M_{\mathrm{act}}-2c_\lambda\nu_\lambda\wideD.
\end{equation}
Analyticity of $\mathcal Z(t,z(t))$ proves convergence of the entire
matrix series, in addition to this two-term expansion.

\emph{Showing that the fast coefficient is nonzero.}
The decisive comparison uses the fixed matrix $\wideC$.
Set
$\Omega_L=\cU(0)^T\cU'(0)$ and
$\Omega_R=\cV(0)^T\cV'(0)$; these are skew-symmetric matrices of
sizes $m\times m$ and $n\times n$.
Subscript $aa$ below denotes their $p\times p$ diagonal blocks
on the active coordinates.
The polar derivative
$\Xi=(\msign A_{\mathrm{act}})'(0)\in\R^{p\times p}$ is skew-symmetric:
differentiate the orthogonality of $\msign A_{\mathrm{act}}(s)$ at $A_{\mathrm{act}}(0)=\Sigma_C$,
where its polar factor is $I_p$.
Differentiating the polar decomposition
$A_{\mathrm{act}}(s)=\msign(A_{\mathrm{act}}(s))(A_{\mathrm{act}}(s)^TA_{\mathrm{act}}(s))^{1/2}$ also gives
$\Sigma_C\Xi+\Xi \Sigma_C=A_{\mathrm{act}}'(0)-A_{\mathrm{act}}'(0)^T$.

In the fixed coordinates $\cU(0),\cV(0)$, the active block of
$\Phi_{\mathrm{part}}'(0)$ is
$-(\Omega_L)_{aa}-\Xi+(\Omega_R)_{aa}$.
This follows by differentiating its three matrix factors:
the active diagonal block of the middle factor at zero is $-I_p$,
and cross-cluster factors have zero active diagonal contribution.
In these same coordinates, $\wideC$ has only the active block $\Sigma_C$.
It follows that
\begin{equation}\label{eq:fixed_C_polar_pairing}
 \begin{aligned}
 \ip{\wideC}{\Phi_{\mathrm{part}}'(0)}
 &=\tr\!\left(
 \Sigma_C[-(\Omega_L)_{aa}-\Xi+(\Omega_R)_{aa}]\right)=0,\\
 \ip{\wideC}{\Phi_{\mathrm{part}}(s)-\bar\Phi}&=O(s^2).
 \end{aligned}
\end{equation}
Here each trace is zero because
$\tr(\Sigma_C T)=\tr((\Sigma_C T)^T)=-\tr(T\Sigma_C)$ for symmetric $\Sigma_C$
and skew-symmetric $T$. This explicitly accounts for the frame
motion as well as the polar derivative.

For later use with general profiles, let $U_a(s)\in\R^{m\times p}$
and $V_a(s)\in\R^{n\times p}$ be the active column groups of
$\cU(s)$ and $\cV(s)$.
For each $s$, choose one full singular value decomposition of $A_{\mathrm{act}}(s)$,
with orthonormal left columns $u_{\mathrm{act},j}(s)\in\R^p$ and right columns
$v_{\mathrm{act},j}(s)\in\R^p$, $1\le j\le p$, labeled so that
$A_{\mathrm{act}}(s)=\sum_{j=1}^p
 k_j(s)u_{\mathrm{act},j}(s)v_{\mathrm{act},j}(s)^T$.
All pairs come from this same orthogonal system, including within
repeated singular subspaces. Define the active matrix outer product
$P_j(s)=(U_a(s)u_{\mathrm{act},j}(s))
 (V_a(s)v_{\mathrm{act},j}(s))^T$.
Since $U_a(s)^T\wideC V_a(s)=\Sigma_C+O(|s|)$ and
$A_{\mathrm{act}}(s)=\Sigma_C+O(|s|)$, uniformly for these unit vectors,
\begin{equation}\label{eq:active_dyad_pairing}
 \begin{aligned}
 \ip{\wideC}{P_j(s)}
 &=u_{\mathrm{act},j}(s)^T[\Sigma_C+O(|s|)]v_{\mathrm{act},j}(s)\\
 &=u_{\mathrm{act},j}(s)^TA_{\mathrm{act}}(s)v_{\mathrm{act},j}(s)+O(|s|)
   =k_j(s)+O(|s|)=\sigma_{C,j}+O(|s|).
 \end{aligned}
\end{equation}
This estimate holds for any choices within repeated singular
subspaces; analyticity of the singular vectors is not required.

Finally, the fixed lifts give
\[
 \ip{\wideC}{M_{\mathrm{act}}}
     =\tr(\Sigma_C\Sigma_C^{-2})=2b_{\ph},\qquad
 \ip{\wideC}{\wideD}=0 .
\]
Pairing \eqref{eq:M} with $\wideC$ now yields
$\ip{\wideC}{M_\Phi}=b_{\ph}/c_\lambda^2=\kv c_\lambda>0$,
which proves \eqref{eq:Cpair}.
If $D\ne0$, the first coefficient is nonzero because
$\nf{\wideD}=\nf{D}$.
If $D=0$, Cauchy--Schwarz gives the nonzero lower bound
$\nf{M_\Phi}\ge\kv c_\lambda/\nf{C}$.
These coefficients prove the two exact direction rates in
\eqref{eq:classification}.
\end{proof}

\subsection{A cubic residual with zero quadratic block}\label{app:cubic_residual}
The $4\times4$ line
\begin{equation}\label{eq:cubicexample}
 A^\star=\diag(1,1,0,0),\qquad
 \Theta=\begin{pmatrix}1&0&0&1\\0&0&0&-1\\0&0&1&0\\1&1&0&0\end{pmatrix}
\end{equation}
has $\Sigma=I_2$ and the blocks
\[
 \Theta_{++}=\diag(1,0),\qquad
 \Theta_{+0}=\begin{pmatrix}0&1\\0&-1\end{pmatrix},\qquad
 \Theta_{0+}=\begin{pmatrix}0&0\\1&1\end{pmatrix},\qquad
 C=\diag(1,0).
\]
Here $U_C=V_C=(1,0)^T$. In the residual coordinate,
$D_0=-(1,1)(1,-1)^T=0$, and hence $D=0$.
The two summands indexed by $i$ in \eqref{eq:kv} are both $2$,
while the skew term from $\Theta_{++}$ is zero, so $\kv=4$.
Its coupled coordinates $1,2,4$ have determinant $s^3$; their two large
singular values tend to one. The remaining singular value is therefore
of order $|s|^3$. Despite this higher mode, the direction rate is
precisely $\delta^{2/3}$ and $\nu_{\lambda}=0$.

\section{Selection, profile laws, and finite-accuracy recovery}
\label{app:profiles}
This appendix proves the results of Section~\ref{sec:profiles}.
We first justify the common limiting direction, then analyze the
tail and saturation balances, and finally account for extrapolation
and bisection errors. All matrices and constants refer to a fixed
input line. Matrix remainder estimates use the Frobenius norm.

\subsection{Why admissible profiles select the same direction}
\label{app:selection}
\begin{proof}[Proof of Proposition~\ref{prop:selection}]
We use scalar Fenchel duality in the form
$\psi(z)+\psi^*(u)\ge zu$, with equality when $u=\psi'(z)$.
Define the finite scalar
$M_\psi=\lim_{z\to\infty}(z-\psi(z))$.
The profile assumptions give $0\le f\le1$ and
\[
 0\le z-\psi(z)=\int_0^z(1-f(u))\,du\le M_\psi
 \quad(z\ge0).
\]
Summing this estimate over the singular values gives
\begin{equation}\label{eq:profile_uniform_bound}
 0\le\nn{A}-\Psi_\delta^\psi(A)\le qM_\psi\delta .
\end{equation}
For the line $A^\star+s\Theta$, the reverse triangle inequality now
bounds its smoothed objective below by
$|s|\nn{\Theta}-\nn{A^\star}-qM_\psi\delta$.
It therefore tends to infinity as $|s|\to\infty$, and a minimizing
$s_\delta^\psi$ exists.

For same-size matrices, von Neumann's singular-value trace
inequality~\cite{mirsky1975} is
\[
 \ip{A}{\Phi}\le\sum_{i=1}^q\sigma_i(A)\sigma_i(\Phi).
\]
Combining it with scalar Fenchel duality gives
\[
 \Psi_\delta^\psi(A)
 =\max_{\Phi\in\R^{m\times n}}
       \{\ip{A}{\Phi}-\delta\Omega_\psi(\Phi)\}.
\]
Equality is attained by aligning the left and right singular vectors
and using the coefficients $f(\sigma_i(A)/\delta)$.
These coefficients vanish at zero, and the even $C^1$ scalar profile
gives a differentiable spectral sum
\cite{lewis1996}. Thus the displayed maximizer is the gradient
$Z_\delta^\psi(A)$ defined in the main text.
At a minimizing $s_\delta^\psi$, the scalar chain rule gives
$\ip{\Theta}{Z_\delta^\psi(A^\star+s_\delta^\psi\Theta)}=0$.
Its matrix gradient therefore belongs to $\cB$.
For each $Q\in\cB$, both terms containing $s_\delta^\psi\Theta$ in
the Fenchel comparison vanish, so this gradient maximizes
\eqref{eq:regularizedprimal}.

To prove its uniqueness, note that $\psi^*$ is even and that
$u\mapsto\psi^*(u)-\mu u^2/2$, extended by infinity outside
$[-1,1]$, is even and convex.
The singular-value convexity theorem \cite{lewis1996} then makes
\[
 \Phi\longmapsto
 \Omega_\psi(\Phi)-\frac{\mu}{2}\nf{\Phi}^2
\]
convex on the matrix space.
Consequently the regularized objective in
\eqref{eq:regularizedprimal} is strictly concave on its convex feasible
set, with at most one maximizer.
Every minimizing scalar multiplier therefore produces the same
$\Phi_\delta^\psi$.

It remains to identify the limit under \eqref{eq:boundary} and
$\kv>0$. Theorem~\ref{thm:rigidity} makes $s=0$ the unique exact
minimizer.
By \eqref{eq:profile_uniform_bound} and minimizing optimality,
$\phi(s_\delta^\psi)\le\phi(0)+qM_\psi\delta$.
Coercivity bounds all these multipliers uniformly for small $\delta$.
Every subsequential limit minimizes $\phi$ and is therefore zero.

For a fixed exact certificate $Q\in\cF$, regularized optimality yields
\begin{equation}\label{eq:profile_selection_gap}
 0\le\phi(0)-\ip{A^\star}{\Phi_\delta^\psi}
 \le\delta\{\Omega_\psi(Q)-\Omega_\psi(\Phi_\delta^\psi)\}.
\end{equation}
Here $\psi^*(1)=M_\psi$: the function $z-\psi(z)$ increases to its
supremum $M_\psi$, and negative $z$ cannot enlarge that supremum.
Evenness gives $\psi^*(-1)=M_\psi$.
The conjugate is closed, convex, and finite on $[-1,1]$, so it is
continuous there, including its finite endpoints.
Thus $\Omega_\psi$ is continuous and bounded on the spectral ball.
Compactness of $\cB$ and \eqref{eq:profile_selection_gap} put every
direction cluster point in $\cF$.
The same inequality also gives
$\Omega_\psi(\Phi_\delta^\psi)\le\Omega_\psi(Q)$, so each cluster
point minimizes $\Omega_\psi$ on $\cF$.

In the boundary-face coordinates of \eqref{eq:faceboundary},
the nonzero singular values consist of $r+p$ ones and the nonzero
singular values of the free residual block $W_{\mathrm{res}}$.
Evenness and strong convexity make zero the unique scalar minimizer
of $\psi^*$, with $\psi^*(0)=0$.
Thus zero singular values contribute nothing, and orthogonal
invariance gives
\[
 \Omega_\psi(\Phi)
 =(r+p)\psi^*(1)
   +\sum_i\psi^*(\sigma_i(W_{\mathrm{res}})),\qquad \Phi\in\cF.
\]
The unique matrix minimum on this face consequently has zero free
residual block and equals $\bar\Phi$.
All direction cluster points coincide, proving the asserted convergence.
\end{proof}

\subsection{The polynomial-tail balance}\label{app:tail}
\begin{proof}[Proof of Theorem~\ref{thm:tail}]
The orthogonal normal form \eqref{eq:normalform} separates three
matrix blocks. We use the positive analytic branches
$k_j(s)$ of $A_{\mathrm{act}}(s)$, $j=1,\ldots,p$, introduced in
Appendix~\ref{app:ph}, with $k_j(0)=\sigma_{C,j}$.
For $s<0$, the active singular values are
$\sigma_{a,j}(s)=-s k_j(s)$.
These are local analytic labels, independent of the global
nonincreasing-order convention.
Their derivatives are
$\sigma_{a,j}'(s)=-k_j(s)-sk_j'(s)$.
The residual block is $R_0(s)=s^2A_{\mathrm{res}}(s)$ in
$\R^{(m_0-p)\times(n_0-p)}$; hence $R_0'(s)=O(|s|)$.

Define the positive scalar coefficient
\[
 c_\psi=
 \left(\frac{\gamma\sum_{j=1}^p \sigma_{C,j}^{1-\rho}}{\kv}\right)^{1/(\rho+1)}.
\]
First test the derivative at $s=-x$, where
$x=z\delta^{\rho/(\rho+1)}$ and $z$ lies in a fixed compact positive interval.
Then $\delta/x\to0$ and $x^2/\delta\to0$.
The active coefficient deficits, meaning their differences from one, are
\[
 d_j(x,\delta):=1-f(\sigma_{a,j}(-x)/\delta)
 =\gamma\delta^\rho x^{-\rho}\sigma_{C,j}^{-\rho}(1+o(1)).
\]
The error is uniform in the chosen interval of $z$ and for the finite
set of $j$.
Since $\sigma_{a,j}'(s)=-\sigma_{C,j}+O(x)$, their correction to the
unsmoothed active trace derivative is
\[
 \sum_{j=1}^p(-\sigma_{a,j}'(-x))d_j
 =\gamma\delta^\rho x^{-\rho}
       \sum_{j=1}^p \sigma_{C,j}^{1-\rho}(1+o(1)).
\]

The large singular values remain bounded below by a fixed positive
number, and their derivatives are bounded.
Their profile deficits therefore contribute $O(\delta^\rho)=o(x)$.
The combined unsmoothed large and active derivative is
$\ell_+'(s)-\ell_{\mathrm{act}}(s)-s\ell_{\mathrm{act}}'(s)=\kv s+O(s^2)$.
For the residual block, \eqref{eq:R2spectralLip} and
$Z_\delta^\psi(0)=0$ give
$\nf{Z_\delta^\psi(R_0(s))}=O(x^2/\delta)$.
Its scalar contribution is the inner product with
$R_0'(s)=O(x)$, and is $O(x^3/\delta)=o(x)$.
Adding the three contributions gives
\begin{equation}\label{eq:tailbalance}
 h_\delta^\psi(\lambda^\star-x)
 =-\kv x+\gamma\delta^\rho x^{-\rho}
                \sum_{j=1}^p \sigma_{C,j}^{1-\rho}+o(x).
\end{equation}
After division by $\delta^{\rho/(\rho+1)}$, the limiting function is
$-\kv z+\gamma z^{-\rho}\sum_j\sigma_{C,j}^{1-\rho}$.
It is strictly decreasing and has its unique positive zero at $c_\psi$.
For fixed $0<z_-<c_\psi<z_+$, the derivatives at
$\lambda^\star-z_+\delta^{\rho/(\rho+1)}$ and
$\lambda^\star-z_-\delta^{\rho/(\rho+1)}$ have opposite signs for small
$\delta$.
Monotonicity brackets every minimizing multiplier between them.
Letting $z_-\uparrow c_\psi$ and $z_+\downarrow c_\psi$ proves
\begin{equation}\label{eq:tail_scalar_coefficient}
 s_\delta^\psi=-c_\psi\delta^{\rho/(\rho+1)}(1+o(1)).
\end{equation}

We next evaluate the matrix direction at any such root.
The local expansion $f(u)=\alpha_0 u+o(u)$, applied to all singular
values of the residual block, gives
\[
 Z_\delta^\psi(s^2A_{\mathrm{res}}(s))
   =\alpha_0\frac{x^2}{\delta}D_0+o(x^2/\delta).
\]
This is a Frobenius-norm estimate: finitely many scalar singular
coefficients obey a common small-argument estimate, including zero
coefficients.
The moving orthogonal frames differ from their values at zero by
$O(x)$, so lifting this term to the original coordinates replaces
$D_0$ by $\wideD$ with a relative $O(x)$ error.
Frame motion in the remaining direction is $O(x)$; active deficits
are $O(\delta^\rho x^{-\rho})=O(x)$; large deficits are $o(x)$.
Because $\delta/x\to0$, all these errors are
$o(x^2/\delta)$.
If $D\ne0$, we conclude that
\[
 \frac{\Phi_\delta^\psi-\bar\Phi}
      {\delta^{(\rho-1)/(\rho+1)}}
 \longrightarrow\alpha_0 c_\psi^2\wideD\ne0.
\]
This proves both the upper and lower bounds of the slow alternative.

If $D=0$, the normal form gives $R_0(s)=O(x^3)$.
Its direction is $O(x^3/\delta)=o(x)$, and the other terms are $O(x)$.
This proves the fast upper bound.
To prove a nonzero lower bound, use the partial polar field
$\Phi_{\mathrm{part}}(s)$ in \eqref{eq:W}.
By \eqref{eq:fixed_C_polar_pairing},
$\ip{\wideC}{\Phi_{\mathrm{part}}(-x)-\bar\Phi}=O(x^2)$.
The active dyads $P_j(s)\in\R^{m\times n}$ defined in
Appendix~\ref{app:ph} satisfy
$\ip{\wideC}{P_j(s)}=\sigma_{C,j}+O(|s|)$ by
\eqref{eq:active_dyad_pairing}.
The active direction differs from the partial polar direction by
$\sum_j d_jP_j(-x)$.
Large deficits and the residual direction are $o(x)$.
Consequently, at a root,
\[
 \ip{\wideC}{\Phi_\delta^\psi-\bar\Phi}
 =\sum_{j=1}^p \sigma_{C,j}d_j+o(x)
 =\gamma\delta^\rho x^{-\rho}\sum_{j=1}^p \sigma_{C,j}^{1-\rho}+o(x)
 =\kv x+o(x),
\]
where the last equality is \eqref{eq:tailbalance}.
After \eqref{eq:tail_scalar_coefficient}, this pairing divided by
$\delta^{\rho/(\rho+1)}$ tends to $\kv c_\psi>0$.
Cauchy--Schwarz then gives the required matching Frobenius lower bound.
\end{proof}

\subsection{Finite saturation and the Huber root}
\label{app:saturation}
\begin{proof}[Proof of Theorem~\ref{thm:saturation}]
\emph{Locating the root.}
In this proof $\tau>0$ is a fixed test ratio of displacement to
$\delta$; it is independent of the cube-root parameter used earlier.
Put $\tau_0=z_\sat/\sigma_{C,\min}$.
At $s=-\tau\delta$ with $0<\tau<\tau_0$, the active scaled singular
values converge to $\tau \sigma_{C,j}$.
The large modes are saturated, and the residual derivative contribution
tends to zero.
Thus
\[
 h_\delta^\psi(\lambda^\star-\tau\delta)
 \longrightarrow\beta-\sum_{j=1}^p \sigma_{C,j} f(\tau \sigma_{C,j})>0 .
\]
At least one mode with $\sigma_{C,j}=\sigma_{C,\min}$ remains below its saturation
threshold, which gives the strict inequality.
For a fixed $\tau>\tau_0$, all large and active modes are saturated
for sufficiently small $\delta$.
Their grouped derivative is
$-\kv\tau\delta+O(\delta^2)$.
The residual block is $O(\delta^2)$ and its derivative is $O(\delta)$,
so the gradient bound \eqref{eq:R2spectralLip} gives a residual scalar
contribution $O(\delta^2)$.
Hence the derivative is negative at this second test point.
Monotonicity, followed by test ratios approaching $\tau_0$ from
both sides, places every minimizing root at
\begin{equation}\label{eq:saturation_localization}
 x:=-s_\delta^\psi=\tau_0\delta+o(\delta).
\end{equation}

\emph{Bounding all active deficits.}
At a root define, for $j=1,\ldots,p$,
\[
 d_j=1-f(xk_j(-x)/\delta)\ge0,\qquad
 d_{\max}=\max_{1\le j\le p}d_j .
\]
The localization makes $d_{\max}\to0$, since
$\tau_0\sigma_{C,j}\ge z_\sat$ and $f$ is continuous.
The correction of active mode $j$ to the partial polar derivative is
$(k_j(s)+sk_j'(s))d_j=(\sigma_{C,j}+O(x))d_j$.
Large modes are saturated.
The residual scalar contribution is
$O(x^3/\delta)=O(x^2)$, because $x$ is comparable to $\delta$.
The root equation therefore becomes
\begin{equation}\label{eq:deficitbalance}
 \sum_{j=1}^p \sigma_{C,j}d_j=\kv x+O(x^2+xd_{\max}).
\end{equation}
Choose a fixed $C_{\rm rem}\ge0$ that bounds the remainder in
absolute value by $C_{\rm rem}(x^2+xd_{\max})$.
The nonnegative sum on the left is at least $\sigma_{C,\min}d_{\max}$, so
\[
 (\sigma_{C,\min}-C_{\rm rem}x)d_{\max}
 \le\kv x+C_{\rm rem}x^2.
\]
For $C_{\rm rem}x\le \sigma_{C,\min}/2$, this gives
$d_{\max}=O(x)$.
Substitution into \eqref{eq:deficitbalance} improves the sum estimate to
$\sum_j\sigma_{C,j}d_j=\kv x+O(x^2)$.

\emph{Upper and lower bounds for the direction.}
Frame motion is $O(x)$, the active deficits are $O(x)$, and the
residual direction is $O(x^2/\delta)=O(\delta)$.
Thus $\nf{\Phi_\delta^\psi-\bar\Phi}=O(\delta)$.
For a lower bound, the partial polar field has fixed-$\wideC$
pairing $O(x^2)$ by \eqref{eq:fixed_C_polar_pairing}.
The active contribution is
$\sum_jd_j\ip{\wideC}{P_j(-x)}
 =\sum_j\sigma_{C,j}d_j+O(x^2)$, using \eqref{eq:active_dyad_pairing}.
At zero, the fixed $\wideC$ has no compression onto the residual
block. Since the frames move by $O(x)$, the residual direction
has pairing $O(x^3/\delta)=O(\delta^2)$.
It follows that
\[
 \ip{\wideC}{\Phi_\delta^\psi-\bar\Phi}
 =\sum_j\sigma_{C,j}d_j+O(\delta^2)
 =\kv x+O(\delta^2).
\]
Together with \eqref{eq:saturation_localization}, this proves
\eqref{eq:satpair}.
Its positive limit and Cauchy--Schwarz yield the matching linear
lower bound in \eqref{eq:satmain}.

\emph{The Huber displacement.}
At a Huber root, at least one active mode is strictly unsaturated.
If all were saturated, every $d_j$ would be zero, and the root
derivative would instead be $-\kv x+O(\delta^2)<0$.
For an unsaturated mode $j$,
$d_j=1-xk_j(-x)/\delta=O(\delta)$, so
\[
 xk_j(-x)=\delta+O(\delta^2),\qquad
 x=\delta/\sigma_{C,j}+O(\delta^2).
\]
Since $x/\delta\to1/\sigma_{C,\min}$ and there are finitely many $\sigma_{C,j}$,
such a mode must have $\sigma_{C,j}=\sigma_{C,\min}$ when $\delta$ is sufficiently small.
This proves $s_\delta^{\rm hub}=-\delta/\sigma_{C,\min}+O(\delta^2)$.

\emph{Uniqueness, including repeated and zero singular values.}
To differentiate the entire spectral sum, choose local singular-value
branches $\widetilde\sigma_1(s),\ldots,\widetilde\sigma_q(s)$ on
a sufficiently small negative interval.
Their labels need not be in decreasing order.
There are $r$ large branches, $p$ active branches
$\sigma_{a,j}(s)=-sk_j(s)$, and $q-r-p$ residual branches.
An analytic signed SVD of the matrix line
\cite{bunsegerstner1991} provides analytic scalar branches.
Each nonzero branch has a fixed sign on a sufficiently small
negative interval; changing that sign once gives a nonnegative
analytic branch there.
Identically zero branches remain zero.
One may also obtain these branches from the smaller analytic Gram
matrix: a nonzero nonnegative analytic eigenvalue has even vanishing
order at zero, so its square root on $s<0$ has an analytic continuation
to zero. This explains why the first two derivatives of all chosen
singular branches are bounded up to zero on that side, independently
of $\delta$.

Fix a sufficiently small $\delta$ and a root.
The strictly unsaturated active mode found above stays below its
clipping edge in an open neighborhood of the root.
Its derivative has magnitude at least $\sigma_{C,\min}/2$, since
$\sigma_{a,j}'(s)=-\sigma_{C,j}+O(|s|)$.
At a point in this neighborhood define the index set
\[
 I_{\rm u}(s,\delta)
 =\{i\in\{1,\ldots,q\}:\widetilde\sigma_i(s)<\delta\}.
\]
On each open interval with fixed clipping status, differentiating the
scalar spectral sum twice gives
\[
 \begin{aligned}
 \frac{d^2}{ds^2}\Psi_\delta^{\rm hub}(A(s))
 &=\sum_{i=1}^q \widetilde\sigma_i''(s)
           f_{\rm hub}(\widetilde\sigma_i(s)/\delta)
   +\frac1\delta\sum_{i\in I_{\rm u}(s,\delta)}
           (\widetilde\sigma_i'(s))^2\\
 &\ge -C_0+\frac{\sigma_{C,\min}^2}{4\delta}>0 .
 \end{aligned}
\]
Here $C_0$ can be chosen as the sum of fixed bounds for
$|\widetilde\sigma_i''|$, since $0\le f_{\rm hub}\le1$.
In a neighborhood where $|s|$ is comparable to $\delta$, all large
branches are saturated and all residual branches are unsaturated.
Each active branch is strictly monotone and therefore has at most
one clipping join.
There are finitely many joins, and the objective is $C^1$ across
them because Huber's scalar profile is $C^1$.
The positive second-derivative bound on the intervening intervals
makes its first derivative strictly increasing through this neighborhood.
The root is consequently a strict local minimum.
A second global minimum of this convex objective would force a flat
segment through the first, a contradiction.
This proves uniqueness.
\end{proof}

\paragraph{The running example}
The explicit Huber expansion for \eqref{eq:slowexample} is
\begin{equation}\label{eq:hubexplicit}
 s_\delta^{\rm hub}=-\delta+2\delta^2+O(\delta^3),\qquad
 \Phi_\delta^{\rm hub}-\bar\Phi
 =\delta\begin{pmatrix}0&0&-1\\0&2&0\\-1&0&-1\end{pmatrix}
       +O(\delta^2).
\end{equation}
In \eqref{eq:slowexample}, coordinates one and three form the
symmetric block
$\left[\begin{smallmatrix}1+s&s\\s&0\end{smallmatrix}\right]$;
coordinate two has eigenvalue $s$.
Call the two block eigenvalues $\vartheta_+(s)$ and
$\vartheta_-(s)$, to distinguish them from scalar multipliers.
Then
$\vartheta_+=1+s+O(s^2)$ and
$\vartheta_-=-s^2+O(s^3)$.
Let $P_+(s)$ denote the orthogonal projector onto the eigenspace of
$\vartheta_+(s)$. Its off-diagonal entry is $s+O(s^2)$ and its lower
diagonal entry is $O(s^2)$.
The theorem just proved places the Huber root at
$s=-\delta+O(\delta^2)$.
Here $\vartheta_+$ is saturated, $\vartheta_-$ is unsaturated,
and the only active mode, $s$, is strictly unsaturated.
For a symmetric input, the singular-value gradient applies the odd
coefficient
$\operatorname{sign}(\vartheta)f_{\rm hub}(|\vartheta|/\delta)$
to its eigenvalues.
Thus the coupled gradient is
$P_+(s)+(\vartheta_-(s)/\delta)(I-P_+(s))$, and the decoupled
coefficient is $s/\delta$.
The tangency equation is
$1+s/\delta+2s+O(\delta^2)=0$, yielding
$s=-\delta+2\delta^2+O(\delta^3)$.
Substituting into those projector and coefficient formulas gives the
matrix in \eqref{eq:hubexplicit}, whose Frobenius norm is $\sqrt7$.

\subsection{A secondary law at the saturation edge}\label{app:edge}
This refinement describes how the scalar multiplier approaches the
threshold in Theorem~\ref{thm:saturation}; the direction remains
linear in $\delta$.

\begin{corollary}[Secondary saturation-edge law]\label{cor:edge}
Under Theorem~\ref{thm:saturation}, let $y=z_\sat-z>0$ denote distance
below the saturation edge, and assume
\begin{equation}\label{eq:edgehyp}
 1-f(z_\sat-y)=\gamma_{\rm e}y^{\theta_{\rm e}}
       +o(y^{\theta_{\rm e}})\quad(y\downarrow0),
 \qquad \gamma_{\rm e}>0,\quad\theta_{\rm e}>1 .
\end{equation}
Let $N=\#\{j\in\{1,\ldots,p\}:\sigma_{C,j}=\sigma_{C,\min}\}$.
Then every minimizing multiplier satisfies
\begin{equation}\label{eq:edgeexp}
 s_\delta^\psi=-\frac{z_\sat}{\sigma_{C,\min}}\delta+
 \frac1{\sigma_{C,\min}}
 \left(\frac{\kv z_\sat}
 {N\gamma_{\rm e}\sigma_{C,\min}^2}\right)^{1/\theta_{\rm e}}
 \delta^{1+1/\theta_{\rm e}}(1+o(1)).
\end{equation}
\end{corollary}
\begin{proof}
Put $x_0=z_\sat\delta/\sigma_{C,\min}$.
For a fixed positive test coefficient $v$, evaluate the derivative at
$s=-x$ with
$x=x_0-v\delta^{1+1/\theta_{\rm e}}$.
For small $\delta$ this is a positive displacement.
For each active mode with $\sigma_{C,j}=\sigma_{C,\min}$, the distance of its scaled
active singular value below the edge is
\[
 y_j=z_\sat-\sigma_{a,j}(-x)/\delta
     =\sigma_{C,\min}v\delta^{1/\theta_{\rm e}}+O(\delta).
\]
Since $\theta_{\rm e}>1$, the remainder is
$o(\delta^{1/\theta_{\rm e}})$.
The corresponding deficit is consequently
$\gamma_{\rm e}\sigma_{C,\min}^{\theta_{\rm e}}v^{\theta_{\rm e}}
 \delta(1+o(1))$, uniformly for $v$ in a fixed compact positive
interval.
All stronger active modes are saturated.
The same grouped derivative calculation used in
\eqref{eq:deficitbalance}, now at these test points, gives
\[
 \frac{h_\delta^\psi(\lambda^\star-x)}{\delta}
 \longrightarrow
 -\frac{\kv z_\sat}{\sigma_{C,\min}}
   +N\gamma_{\rm e}\sigma_{C,\min}^{\theta_{\rm e}+1}v^{\theta_{\rm e}} .
\]
This function has the unique positive zero
$v_{\rm e}=(\kv z_\sat/
 (N\gamma_{\rm e}\sigma_{C,\min}^{\theta_{\rm e}+2}))^{1/\theta_{\rm e}}$.
Fixed $0<v_-<v_{\rm e}<v_+$ give opposite derivative signs;
the corresponding $s$-values are
$-x_0+v_\pm\delta^{1+1/\theta_{\rm e}}$ in increasing order.
Monotonicity brackets every root between them.
Letting $v_\pm$ approach $v_{\rm e}$ proves
\eqref{eq:edgeexp}.
\end{proof}

For example, the odd derivative obtained from
$f(z)=2z-z^2$ on $[0,1]$ and $f(z)=1$ for $z\ge1$ is
2-Lipschitz.
Its even integral profile has a $1/2$-strongly convex conjugate.
Here $z_\sat=\gamma_{\rm e}=1$ and $\theta_{\rm e}=2$, so the
secondary displacement has order $\delta^{3/2}$.
At exponent one, the motion of $A_{\mathrm{act}}(s)$ also affects the
order-$\delta^2$ coefficient, as in the Huber calculation above.

\subsection{Extrapolation and radial repair}\label{app:extrapolation}
\begin{proof}[Proof of Corollary~\ref{cor:extrap}]
By Theorem~\ref{thm:PH},
$\Phi_\delta=\bar\Phi+\sum_{\ell\ge1}M_\ell t^\ell$,
where $t=\delta^{1/3}$ and each $M_\ell\in\R^{m\times n}$ is fixed.
Its first two coefficients are $M_1=c_\lambda^2\wideD$ and $M_2=M_\Phi$.
At scales $\delta/8$ and $\delta/64$ the series parameters are $t/2$
and $t/4$.
The weights $(-1,2)$ sum to one and cancel the linear coefficient.
For the three-scale weights, direct calculation gives
\[
 (1-6+8)/3=1,\qquad
 (1-6/2+8/4)/3=0,\qquad
 (1-6/4+8/16)/3=0 .
\]
The raw errors are therefore $O(t^2)$ and $O(t^3)$.
Each raw combination is tangent because every exact input direction
is tangent.

For $Q\in\R^{m\times n}$, put $\xi=\max(1,\no{Q})$.
Since $\no{\bar\Phi}\le1$,
$0\le\xi-1\le\no{Q-\bar\Phi}\le\nf{Q-\bar\Phi}$.
The identity
\[
 Q/\xi-\bar\Phi=(Q-\bar\Phi)/\xi-(1-1/\xi)\bar\Phi
\]
gives \eqref{eq:repairbound}.
Radial scaling preserves tangency and enforces the spectral bound.
It therefore preserves both raw error orders.

More generally, fix an integer $k\ge1$, distinct positive scale factors
$\zeta_0,\ldots,\zeta_k$, and real weights $w_0,\ldots,w_k$ satisfying
$\sum_{j=0}^kw_j=1$ and
$\sum_{j=0}^kw_j\zeta_j^{\ell/3}=0$ for $\ell=1,\ldots,k$.
Then $\mathcal N(\sum_{j=0}^kw_j\Phi_{\zeta_j\delta})$ has error
$O(\delta^{(k+1)/3})$.

The positive nodes $\zeta_j^{1/3}$ are distinct.
The $k+1$ moment conditions, including the sum-one condition, form
a square Vandermonde system for $w_0,\ldots,w_k$.
They cancel the first $k$ powers of $t$ and leave
$O(t^{k+1})=O(\delta^{(k+1)/3})$.
The same repair inequality proves the general assertion.
\end{proof}

\subsection{The bisection direction solver}\label{app:oracle}
The strong-convexity step underlying bisection is useful for all
admissible profiles.
For two same-size matrices $A_1,A_2$, write
$Q_i=Z_\delta^\psi(A_i)$, $i=1,2$.
The matrix Fenchel maximizers obey
$A_i/\delta\in\partial\Omega_\psi(Q_i)$.
For a $\mu$-strongly convex function, subgradients at two points
satisfy
\[
 \ip{A_1/\delta-A_2/\delta}{Q_1-Q_2}
 \ge\mu\nf{Q_1-Q_2}^2 .
\]
This follows by adding the two strong subgradient inequalities,
one with base point $Q_1$ and one with base point $Q_2$.
Cauchy--Schwarz and division by $\nf{Q_1-Q_2}$ when it is nonzero
give \eqref{eq:R2spectralLip}; the zero case is immediate.
The extension of $\Omega_\psi$ by infinity incorporates the spectral
ball constraint in this subgradient formula.
For pseudo-Huber the scalar conjugate is
$\psi_{\ph}^*(u)=1-\sqrt{1-u^2}$ on $[-1,1]$, and for Huber it is $u^2/2$
there. Each is 1-strongly convex, which explains the common
$\mu=1$ used by the algorithm.

\begin{proof}[Proof of Proposition~\ref{prop:R2oracle}]
On the input line, let
$\lambda_\delta^\psi=\lambda^\star+s_\delta^\psi$ be any minimizing
multiplier.
The objective is differentiable, convex, and coercive, as shown in
Appendix~\ref{app:selection}.
Its derivative $h_\delta^\psi$ is continuous and nondecreasing.
It takes a negative and a positive value: if it were nonnegative
everywhere, the objective on the negative half-line would be bounded
above by its value at zero, contradicting coercivity.
The everywhere nonpositive case contradicts coercivity on the
positive half-line.
Symmetric bracket doubling therefore terminates.

Every bisection preserves a minimizing root and halves the interval.
A zero derivative gives the exact matrix direction by
Proposition~\ref{prop:selection}.
Otherwise, if the final bracket width is $w$, its midpoint is within
$w/2$ of a minimizing root.
The gradient bound with $\mu=1$ gives
\[
 \nf{\widetilde\Phi_\delta^\psi-\Phi_\delta^\psi}
 \le\frac{\nf{\Theta}}{\delta}\frac w2\le\epsin .
\]
Starting from width $w_{\rm br}$, solving
$2^{-j}w_{\rm br}\le2\delta\epsin/\nf{\Theta}$
gives the stated bisection count.
Tangency of the exact $\Phi_\delta^\psi$ gives
$|\ip{\Theta}{\widetilde\Phi_\delta^\psi}|
 \le\nf{\Theta}\epsin$.
All gradient singular coefficients belong to $[0,1]$, so the
computed direction has operator norm at most one.
\end{proof}

\subsection{Combining the bias, bisection error, and feasibility repair}
\label{app:recovery}
\begin{proof}[Proof of Proposition~\ref{prop:R2recovery}]
The centered/input coordinate identity
\eqref{eq:R3coordinates} identifies the input line with the exact
paths of the preceding theorems.
Let $Q^{\rm ex}$ be the exact Huber direction in mode H or the
exact three-scale raw combination in mode E.
Theorem~\ref{thm:saturation} or
Corollary~\ref{cor:extrap} respectively gives
$\nf{Q^{\rm ex}-\bar\Phi}\le C_{\rm rec}\delta$.
The constant depends only on the fixed line and chosen mode, on
a sufficiently small fixed interval of $\delta$.

For mode H the bisection error is at most $\epsin$.
For mode E, the absolute combination weights sum to
$(1+6+8)/3=5$.
The triangle inequality therefore gives
\[
 \nf{Q-Q^{\rm ex}}\le\omega_{\rm mode}\epsin,\qquad
 \nf{Q-\bar\Phi}\le C_{\rm rec}\delta+
                         \omega_{\rm mode}\epsin .
\]
The map $\Pi_\Theta$ is the orthogonal projection in the Frobenius
matrix space onto $\{Q:\ip{\Theta}{Q}=0\}$, and is nonexpansive.
Since it fixes $\bar\Phi$,
$\nf{\Pi_\Theta(Q)-\bar\Phi}\le\nf{Q-\bar\Phi}$.
Applying \eqref{eq:repairbound}, with
$\nf{\bar\Phi}=\sqrt{r+p}\le\sqrt q$, proves
\eqref{eq:R2finiteaccuracy}.
Projection and radial normalization give the two exact feasibility
conditions directly.

For stability of repair, set
$Y=\Pi_\Theta(Q)$, $Y'=\Pi_\Theta(Q')$,
$\xi=\max(1,\no{Y})$, and $\xi'=\max(1,\no{Y'})$.
Both denominators are at least one,
$|\xi-\xi'|\le\nf{Y-Y'}$, and
$\nf{Y'}/\xi'\le\sqrt q$.
Therefore
\[
 \begin{aligned}
 \nf{Y/\xi-Y'/\xi'}
 &\le\frac{\nf{Y-Y'}}{\xi}
       +\frac{\nf{Y'}}{\xi'}\frac{|\xi-\xi'|}{\xi}\\
 &\le(1+\sqrt q)\nf{Y-Y'}
 \le(1+\sqrt q)\nf{Q-Q'} .
 \end{aligned}
\]
This is \eqref{eq:R2repairLip}.
\end{proof}

\section{Spectral tools and the initial selection proof}\label{app:foundations}
This appendix completes Proposition~\ref{prop:PHselection}.
It supplies the precise spectral Hessian and convex lifting facts used
there.

\Needspace{4\baselineskip}
\subsection{Existence, uniqueness, and the selected center}\label{app:ph_selection_proof}
\begin{proof}[Full proof of Proposition~\ref{prop:PHselection}]
The bound $0\le\nn{A}-\Psi_\delta(A)\le q\delta$ gives
\[
 \Psi_\delta(A^\star+s\Theta)
 \ge |s|\nn{\Theta}-\nn{A^\star}-q\delta .
\]
The right-hand side tends to infinity with $|s|$, so a minimum exists.

To show uniqueness we use the symmetric dilation $\mathcal H(A)$
defined in Appendix~\ref{app:normalform} and the scalar function
\[
 g_\delta(u)=\sqrt{u^2+\delta^2}\quad(u\in\R).
\]
The eigenvalues of $\mathcal H(A)$ are the $q$ pairs
$\pm\sigma_i(A)$ and $m+n-2q$ additional zeros.
For a real symmetric matrix, $g_\delta$ acts on its eigenvalues.
Consequently
\begin{equation}\label{eq:PHdilation}
 \Psi_\delta(A)=\tfrac12\tr g_\delta(\mathcal H(A))
                         -\tfrac{m+n}{2}\delta .
\end{equation}
This identity includes rectangular inputs and singular matrices.

We use the trace Hessian in the form needed for this line
\cite{lewis2001}. Choose an orthogonal eigenbasis
$Q_{\mathrm{eig}}\in\R^{(m+n)\times(m+n)}$ for
$\mathcal H(A)$, write its eigenvalues as
$\xi_1,\ldots,\xi_{m+n}$, and set
$\widetilde\Theta
 =Q_{\mathrm{eig}}^T\mathcal H(\Theta)Q_{\mathrm{eig}}$.
Thus $\widetilde\Theta\in\R^{(m+n)\times(m+n)}$ contains the
dilated direction entries in this eigenbasis.
The second derivative is
\[
 \begin{aligned}
 &D^2\{\tfrac12\tr g_\delta\}(\mathcal H(A))
       [\mathcal H(\Theta),\mathcal H(\Theta)]\\
 &\qquad=\tfrac12\sum_{i=1}^{m+n} g_\delta''(\xi_i)\widetilde\Theta_{ii}^{\,2}
  +\sum_{1\le i<j\le m+n}
    \frac{g_\delta'(\xi_i)-g_\delta'(\xi_j)}{\xi_i-\xi_j}
            \widetilde\Theta_{ij}^{\,2}.
 \end{aligned}
\]
At a repeated eigenvalue the quotient means $g_\delta''(\xi_i)$.
Here $D^2$ denotes the second directional derivative; the indices refer
to the transformed matrix $\widetilde\Theta$.
Since $g_\delta''(u)=\delta^2(u^2+\delta^2)^{-3/2}>0$,
every coefficient is positive.
Apply this formula at $A=A^\star+s\Theta$.
Since $\Theta\ne0$, its symmetric direction $\mathcal H(\Theta)$ is nonzero.
At least one entry of $\widetilde\Theta$ is nonzero, so the line objective
has strictly positive second derivative everywhere.
Its unique minimizing point is $s_\delta$.
Differentiation of the line gives
$\ip{\Theta}{Z_\delta(A^\star+s_\delta\Theta)}=0$.
All gradient singular coefficients are less than one, proving tangency
and $\no{\Phi_\delta}<1$.

The Fenchel maximizer in \eqref{eq:PHfenchel} at
$A^\star+s_\delta\Theta$ is $\Phi_\delta$.
For any $Q\in\cB$, tangency cancels the term with $s_\delta\Theta$,
and Fenchel optimality gives
\[
 \ip{A^\star}{\Phi_\delta}-\delta\Omega_{\ph}(\Phi_\delta)
 \ge \ip{A^\star}{Q}-\delta\Omega_{\ph}(Q).
\]
To check uniqueness in this matrix problem, use the scalar pseudo-Huber
conjugate $\psi_{\ph}^*(u)=1-\sqrt{1-u^2}$ on $[-1,1]$,
extended by infinity outside. The scalar function
$\psi_{\ph}^*(u)-u^2/2$ is even and convex:
on $(-1,1)$ its second derivative is
$(1-u^2)^{-3/2}-1\ge0$, and it extends continuously to the endpoints.
The singular-value lift of an even convex scalar function is convex
\cite{lewis1996}. Since $\nf{\Phi}^2=\sum_i\sigma_i(\Phi)^2$,
this yields convexity of
$\Omega_{\ph}(\Phi)-\nf{\Phi}^2/2$.
Equivalently, $\Omega_{\ph}$ is 1-strongly convex in Frobenius norm.
The regularized maximizer on the convex set $\cB$ is therefore unique.

Now assume \eqref{eq:boundary} and a unique exact multiplier.
Approximation and optimality imply
$\phi(s_\delta)\le\phi(0)+q\delta$.
Coercivity bounds the minimizing points for small $\delta$, and every
cluster point minimizes $\phi$; uniqueness gives $s_\delta\to0$.
For a fixed $Q\in\cF$, the same optimality comparison gives
\[
 0\le\ip{A^\star}{Q-\Phi_\delta}
 \le\delta\{\Omega_{\ph}(Q)-\Omega_{\ph}(\Phi_\delta)\}.
\]
The feasible set $\cB$ is compact and $\Omega_{\ph}$ is continuous
there. Thus each cluster point of $\Phi_\delta$ lies in $\cF$
and minimizes $\Omega_{\ph}$ on that face.
For a member of \eqref{eq:faceboundary}, the nonzero singular values
consist of $r+p$ ones and the nonzero singular values of
$W_{\mathrm{res}}$.
The scalar conjugate $\psi_{\ph}^*$ has its unique minimum at zero. Hence
\[
 \Omega_{\ph}(\Phi)
  =(r+p)\psi_{\ph}^*(1)
    +\sum_i\psi_{\ph}^*(\sigma_i(W_{\mathrm{res}}))
\]
is uniquely minimized by the zero residual block.
The selected limit is exactly $\bar\Phi$.
\end{proof}

\section{Inputs and computational protocols}\label{app:protocol}
This appendix specifies the fixed test matrices and gives the settings
and results for the asymptotic and finite-accuracy tests.
Matrices are listed in their construction coordinates. Ordered singular
values, such as $\sigma_1$ in Corollary~\ref{cor:rankone}, are indexed by
decreasing size. When passing to the ordered coordinates of
\eqref{eq:blocks}, apply the corresponding row and column permutations
to $A^\star$ and every affected block of $\Theta$.

\subsection{The ten fixed matrix lines}\label{app:testinputs}
For each case in Table~\ref{tab:testinputs}, let $m,n,r,p$ have
the displayed values and set $m_0=m-r$ and $n_0=n-r$. The positive diagonal
matrix $\Sigma$ has the listed diagonal entries.
The list $(\sigma_{C,1},\ldots,\sigma_{C,p})$ specifies
$C\in\R^{m_0\times n_0}$ by
\[
 C_{ij}=
 \begin{cases}
  \sigma_{C,i},&i=j\le p,\\
  0,&\text{otherwise},
 \end{cases}
 \qquad 1\le i\le m_0,\quad 1\le j\le n_0 .
\]
For the explicitly listed blocks $\Theta_{++},\Theta_{+0},\Theta_{0+}$ below, the full inputs are
\begin{equation}\label{eq:testconstruction}
 A^\star=\begin{pmatrix}\Sigma&0_{r\times n_0}\\
                   0_{m_0\times r}&0_{m_0\times n_0}\end{pmatrix},
 \qquad
 \Theta=\begin{pmatrix}\Theta_{++}&\Theta_{+0}\\\Theta_{0+}&C\end{pmatrix}.
\end{equation}
Thus $\Theta_{++}\in\R^{r\times r}$, $\Theta_{+0}\in\R^{r\times n_0}$ and
$\Theta_{0+}\in\R^{m_0\times r}$.
Every case has $\tr \Theta_{++}=\sum_{j=1}^p \sigma_{C,j}>0$ and $\kv>0$.
For these verification inputs we set $G=A^\star$ and
$\lambda^\star=0$, so $\lambda=s$ in \eqref{eq:R3coordinates}.
In these bases the reference direction is
$\bar\Phi=\diag(I_r,-I_p,0_{(m_0-p)\times(n_0-p)})$,
with rectangular block-diagonal notation and empty blocks allowed.

\begin{table}[htbp]
\centering
\caption{Parameters for \eqref{eq:testconstruction}.
The last column states whether the residual quadratic block is nonzero.}
\label{tab:testinputs}
\begin{tabular}{lcccccc}
\toprule
Case & $m\times n$ & $r$ & $p$ & diagonal of $\Sigma$ &
$(\sigma_{C,1},\ldots,\sigma_{C,p})$ & $D\ne0$\\
\midrule
P01 & $3\times 3$ & $1$ & $1$ & $(1)$ & $(1)$ & yes\\
P02 & $3\times 3$ & $1$ & $1$ & $(1)$ & $(1)$ & yes\\
P03 & $4\times 4$ & $2$ & $1$ & $(1,2)$ & $(1)$ & no\\
P04 & $5\times 5$ & $2$ & $2$ & $(1,2)$ & $(1,1)$ & yes\\
P05 & $6\times 5$ & $2$ & $2$ & $(1,2)$ & $(1,2)$ & yes\\
P06 & $4\times 6$ & $2$ & $1$ & $(1,2)$ & $(1)$ & yes\\
P07 & $5\times 4$ & $2$ & $2$ & $(1,2)$ & $(1,2)$ & no\\
P08 & $4\times 4$ & $2$ & $2$ & $(1,2)$ & $(1,1)$ & no\\
P09 & $4\times 5$ & $1$ & $1$ & $(1)$ & $(1)$ & no\\
P10 & $4\times 4$ & $2$ & $1$ & $(1,1)$ & $(1)$ & no\\
\bottomrule
\end{tabular}
\end{table}
\paragraph{P01: The running quadratic-residual example}
\[
 \Theta_{++}=\begin{pmatrix}1\end{pmatrix},\qquad
 \Theta_{0+}=\begin{pmatrix}0\\1\end{pmatrix},\qquad
 \Theta_{+0}=\begin{pmatrix}0&1\end{pmatrix}.
\]
\paragraph{P02: A nonsymmetric sign change}
\[
 \Theta_{++}=\begin{pmatrix}1\end{pmatrix},\qquad
 \Theta_{0+}=\begin{pmatrix}0\\-1\end{pmatrix},\qquad
 \Theta_{+0}=\begin{pmatrix}0&1\end{pmatrix}.
\]
\paragraph{P03: A zero quadratic block with a cubic mode}
\[
 \Theta_{++}=\begin{pmatrix}1&0\\0&0\end{pmatrix},\qquad
 \Theta_{0+}=\begin{pmatrix}0&0\\1&2\end{pmatrix},\qquad
 \Theta_{+0}=\begin{pmatrix}0&1\\0&-1\end{pmatrix}.
\]
This test uses $\Sigma=\diag(1,2)$ and $\kv=9/2$.
Its residual cancellation is
\[
 -(1,2)\diag(1,1/2)(1,-1)^T=0.
\]
It is the nonuniform-positive-spectrum variant of
\eqref{eq:cubicexample}; that displayed example instead uses
$\Sigma=I_2$ and a different $\Theta_{0+}$ block.
For P03, $\det(A^\star+s\Theta)=2s^4$:
the two large singular values have a nonzero product limit and
one active singular value is of order $|s|$, leaving a cubic mode.
\paragraph{P04: Repeated active singular values}
\[
 \Theta_{++}=\begin{pmatrix}1&2\\-1&1\end{pmatrix},\qquad
 \Theta_{0+}=\begin{pmatrix}1&0\\0&1\\1&-1\end{pmatrix},\qquad
 \Theta_{+0}=\begin{pmatrix}0&1&1\\1&0&1\end{pmatrix}.
\]
\paragraph{P05: A tall matrix with two active modes}
\[
 \Theta_{++}=\begin{pmatrix}2&1\\0&1\end{pmatrix},\qquad
 \Theta_{0+}=\begin{pmatrix}1&0\\0&1\\1&2\\-1&1\end{pmatrix},\qquad
 \Theta_{+0}=\begin{pmatrix}1&0&2\\0&1&-1\end{pmatrix}.
\]
\paragraph{P06: A wide nonsymmetric matrix}
\[
 \Theta_{++}=\begin{pmatrix}1&1\\-1&0\end{pmatrix},\qquad
 \Theta_{0+}=\begin{pmatrix}1&0\\1&2\end{pmatrix},\qquad
 \Theta_{+0}=\begin{pmatrix}1&1&0&2\\0&-1&1&0\end{pmatrix}.
\]
\paragraph{P07: Full column rank in the compressed normal}
\[
 \Theta_{++}=\begin{pmatrix}1&2\\-1&2\end{pmatrix},\qquad
 \Theta_{0+}=\begin{pmatrix}1&0\\0&1\\1&1\end{pmatrix},\qquad
 \Theta_{+0}=\begin{pmatrix}0&1\\1&0\end{pmatrix}.
\]
\paragraph{P08: A full compressed block with repeated singular values}
\[
 \Theta_{++}=\begin{pmatrix}1&1\\-1&1\end{pmatrix},\qquad
 \Theta_{0+}=\begin{pmatrix}1&0\\0&1\end{pmatrix},\qquad
 \Theta_{+0}=\begin{pmatrix}0&1\\1&0\end{pmatrix}.
\]
\paragraph{P09: A rank-one normal with large nullities}
\[
 \Theta_{++}=\begin{pmatrix}1\end{pmatrix},\qquad
 \Theta_{0+}=\begin{pmatrix}1\\0\\0\end{pmatrix},\qquad
 \Theta_{+0}=\begin{pmatrix}1&0&0&0\end{pmatrix}.
\]
\paragraph{P10: A rank-two normal with repeated positive singular values}
\[
 \Theta_{++}=\begin{pmatrix}1&1\\0&0\end{pmatrix},\qquad
 \Theta_{0+}=\begin{pmatrix}1&0\\0&0\end{pmatrix},\qquad
 \Theta_{+0}=\begin{pmatrix}1&0\\0&0\end{pmatrix}.
\]

\subsection{Asymptotic and structural verification}\label{app:asymptotic_tests}
The structured spectra and scalar roots are computed with 90 decimal digits.
For pseudo-Huber, it evaluates every fixed case at
$\delta=10^{-6},10^{-12},10^{-18},10^{-24}$.
The root tests are initialized at $-1.5$ and $-0.5$ times the
predicted positive displacement magnitude; the derivative is checked
to have opposite signs there.
Pseudo-Huber and polynomial-tail roots are computed by a secant solve
with tolerance $10^{-78}$.
Finite-saturation roots use 255 bisection steps in the corresponding
predicted-scale bracket.
Each root is accepted only if its absolute tangency residual is below
$10^{-60}$.
These brackets test the asymptotic prediction; the algorithmic tests
below use the data-based doubling bracket instead.

For the polynomial-tail family $f_\rho$ in
Section~\ref{sec:profiles}, P01, P03 and P06 are evaluated at
$\rho=3/2,4$ and $\delta=10^{-12},10^{-24},10^{-36}$.
Huber tests use P01, P03, P05, P06 and P08 at
$\delta=10^{-5},10^{-9}$.
The odd derivative obtained from $f(z)=2z-z^2$ on $[0,1]$ and
$f(z)=1$ above one is tested on P01, P03 and P05 at
$\delta=10^{-6},10^{-12}$.
Here the saturation-edge exponent is $\theta_{\rm e}=2$ in the notation
of Corollary~\ref{cor:edge} in Appendix~\ref{app:edge}.
On P01, at $\delta=10^{-12}$, the scaled displacement correction
$(s_\delta^\psi+\delta)/\delta^{3/2}$ is $1.414212562$,
approaching the predicted $\sqrt2$. Appendix~\ref{app:closed_form_check}
checks the same correction independently on P01.

Two- and three-scale pseudo-Huber extrapolations on P01 use
base scales $\delta=10^{-6},10^{-12},10^{-18}$ and the weights in
\eqref{eq:extrapdef}.
The tests measure both the raw spectral norms and the errors
after radial repair. All comparisons use the explicit
$\bar\Phi$ from the input construction.

The matrix-coefficient tests also evaluate the first two terms
in \eqref{eq:Pexp}.
The derivative of the partial polar field $\Phi_{\mathrm{part}}$ is
approximated by the one-sided three-point difference
\[
 \Phi_{\mathrm{part}}'(0)\approx
 \frac{3\bar\Phi-4\Phi_{\mathrm{part}}(-h_{\rm fd})+
              \Phi_{\mathrm{part}}(-2h_{\rm fd})}{2h_{\rm fd}},
 \qquad h_{\rm fd}=10^{-20}.
\]
Here $h_{\rm fd}$ is the finite-difference step.
This estimate is substituted into \eqref{eq:M} and compared with
the scaled second matrix coefficient.

The separate double-precision structural checks use random seed $10092026$.
They construct 120 integer-block boundary examples, of which 119
have $\kv>0$ and one is flat.
For each isolated case they check the curvature and rank inequalities,
the multiplier scale, and orthogonal equivariance.
The five sharp rank examples use the construction of
Theorem~\ref{thm:rank}; Table~\ref{tab:rankchecks} lists their ranks and curvatures.
The fifteen additional rank-one curvature tests have
$A^\star=\diag(1,2,0,0)$ and $\Theta=u_\ell v^T$, where
$u_\ell=(1,\ell,1,0)^T\in\R^4$,
$v=(1,0,1,0)^T\in\R^4$, and $\ell=1,\ldots,15$.
Their $r=2$ positive block is fixed while the factor varies.

The flat control has
\[
 A^\star=\diag(1,0,0),\qquad
 \Theta=\begin{pmatrix}1&1&0\\-1&1&0\\0&0&0\end{pmatrix},
\]
and $\kv=D=0$.
Its negative branch is constant near zero.
The zero-boundary control instead takes
$A^\star=\diag(1,0)$ and
$\Theta=\left[\begin{smallmatrix}0&1\\1&0\end{smallmatrix}\right]$.
It has $a=\beta=0$ and a pseudo-Huber root at $s=0$ by symmetry.

\paragraph{Recorded coefficients and recovery errors}
Between $\delta_1=10^{-18}$ and $\delta_2=10^{-24}$, the five
$D\ne0$ cases give exponents from $0.33333331$ to $0.33333333$;
the five $D=0$ cases give $0.666666667$ to nine decimal places.
These are evaluated by \eqref{eq:observedslope}.

For the running example \eqref{eq:slowexample}, at $\delta=10^{-24}$,
the normalized displacement and its next correction are
$-s_\delta/\delta^{2/3}=0.629960524$ and
$(s_\delta+c_\lambda\delta^{2/3})/\delta=0.083333339$.
They match $c_\lambda=4^{-1/3}$ and $\nu_\lambda=1/12$ in
\eqref{eq:slowconstants}.
The coefficient check
$\nf{(\Phi_\delta-\bar\Phi)/\delta^{1/3}-c_\lambda^2\wideD}
 =1.55\times10^{-8}$
also tests the matrix leading term, while
$\ip{\wideC}{\Phi_\delta-\bar\Phi}/\delta^{2/3}$
approaches $\kv c_\lambda=1.259921050$.
Thus the tests distinguish the partial curvature $\kv=2$ from
the full contact curvature $\kf=4$, which would give the displacement
constant $1/2$.

On P01, tail exponents $\rho=3/2$ and $4$ give observed direction
exponents $0.19999981$ and $0.60000000$, matching $1/5$ and $3/5$
in Theorem~\ref{thm:tail}.
Huber gives $0.99999969$, as predicted by
Theorem~\ref{thm:saturation}.
Every error uses the same explicit $\bar\Phi=\diag(1,-1,0)$.
For extrapolation, at $\delta=10^{-6}$ the raw two-scale combination
has spectral norm $1.000063245$.
Radial repair restores feasibility, giving error
$7.76523\times10^{-5}$; the repaired three-scale error is
$7.77328\times10^{-8}$.
These values check Corollary~\ref{cor:extrap}, including the need to
repair a combination with negative weights.

\begin{table}[htbp]
\centering
\caption{Equality cases for the sharp rank bound. Here $p=\rank C$,
$d=\rank D$, and the normal is $\Theta$.}
\label{tab:rankchecks}
\begin{tabular}{rrrr}
\toprule
$p$ & $d$ & $\rank\Theta=p+2d$ & $\kv=2d$\\
\midrule
1&1&3&2\\
2&1&4&2\\
1&2&5&4\\
3&2&7&4\\
2&3&8&6\\
\bottomrule
\end{tabular}
\end{table}

\subsection{Independent evaluation of the running example}\label{app:closed_form_check}
The independent calculation uses 120 decimal digits and evaluates
the odd spectral coefficient on the two explicit eigenvalues of
the coupled coordinates $1,3$ of P01 and the decoupled eigenvalue $s$.
A 400-step monotone bisection, stopped earlier if adjacent
working-precision numbers coincide, starts from $[-1/4,0]$.
At the same four scales $10^{-6},10^{-12},10^{-18},10^{-24}$,
the calculation checks the pseudo-Huber first and second displacement
coefficients, the Huber linear and quadratic displacement
coefficients, and the rounded-clipping edge correction.
Each of its 12 roots has tangency residual below $10^{-80}$.
The pseudo-Huber results are compared with the singular-value-based
calculation in the preceding subsection.

\subsection{Finite-accuracy recovery and update checks}\label{app:update_checks}
The tests use the two explicit inputs of Section~\ref{sec:numapp}
and 90 decimal digits.
For each $\delta=10^{-3},10^{-6},10^{-9}$ it evaluates three
constructions: one pseudo-Huber direction, one Huber direction,
and the three-scale pseudo-Huber combination.
Every scalar solve uses Algorithm~\ref{alg:oracle}, starting at
bracket radius one and doubling until the endpoint derivatives
bracket zero.
It stops at interval width
$2\delta\epsin/\nf{\Theta}$, with $\delta$ replaced by the actual
smoothing scale in each multiscale call.
The fine solve has $\epsin=\delta^3$ and the corresponding coarse
solve has $\epsin=\delta/100$.
All constructions apply \eqref{eq:R2repair} before any error or
objective comparison.

For each final output $\widehat\Phi$, the recorded quantities are
\[
 \nf{\widehat\Phi-\bar\Phi},\quad
 |\ip{\Theta}{\widehat\Phi}|,\quad
 \no{\widehat\Phi},\quad
 v^\star-\ip{G}{\widehat\Phi},\quad
 \mathcal L(X_0+\eta\widehat\Phi),
 \qquad \eta=0.1 .
\]
Here $G=Y$ only at the chosen start $X_0=0$, and
$\mathcal L(X)=\tfrac12\nf{X-Y}^2$.
The exact identity
\[
 \mathcal L(\eta\widehat\Phi)=
 \tfrac12\nf{Y}^2-\eta\ip{Y}{\widehat\Phi}
                  +\tfrac12\eta^2\nf{\widehat\Phi}^2
\]
checks the update independently of the descent bound.
Its largest recorded discrepancy across the 18 fine runs is below
$10^{-70}$.
The tests also check the direction-value bound,
\eqref{eq:R2decrease}, and the coarse--fine repair bound
\eqref{eq:R2repairLip} using weight sums one and five.
The rank-one three-scale observed exponent is approximately
$1.33330$; this additional cancellation is consistent with the
general upper error bound in Corollary~\ref{cor:extrap}.

\begin{table}[htbp]
\centering
\caption{Final errors $e_{\rm dir}=\nf{\widehat\Phi-\bar\Phi}$ after tangent
projection and radial repair, with $\epsin=\delta^3$.
Huber and three-scale recovery test Proposition~\ref{prop:R2recovery}.}
\label{tab:R2application}
\begin{tabular}{llrrr}
\toprule
$\rank\Theta$ & $\delta$ & Pseudo-Huber & Huber & Three-scale\\
\midrule
3 & $10^{-3}$ & $4.1280\times10^{-2}$ & $2.6450\times10^{-3}$ & $7.7886\times10^{-5}$\\
3 & $10^{-9}$ & $3.9675\times10^{-4}$ & $2.6458\times10^{-9}$ & $7.7704\times10^{-11}$\\
1 & $10^{-3}$ & $2.4054\times10^{-2}$ & $4.8666\times10^{-3}$ & $1.0438\times10^{-4}$\\
1 & $10^{-9}$ & $2.4495\times10^{-6}$ & $4.8990\times10^{-9}$ & $1.0716\times10^{-12}$\\
\bottomrule
\end{tabular}
\end{table}

\paragraph{The error plot}
Figure~\ref{fig:boundary_recovery} uses the three scales
$\delta=10^{-3},10^{-6},10^{-9}$ from the application tests.
The solid curves join the final direction errors. The dashed segments
are slope guides with powers $1/3$, $2/3$, and $1$.

\end{document}